\documentclass[11pt, a4paper]{article}

\usepackage{amsmath}
\usepackage{amssymb}
\usepackage{amsthm}
\usepackage{txfonts}
\usepackage{stmaryrd}
\usepackage{subfig}
\usepackage{graphicx}
\usepackage{tikz}

\usepackage[auth-sc]{authblk}

\usepackage[utf8]{inputenc}
\usepackage[english]{babel}
\usepackage[T1]{fontenc}
\usepackage{lmodern}
\usepackage{microtype}
\usepackage[
    linktocpage=true,
    colorlinks=true,
    linkcolor=purple,
    citecolor=black,
    urlcolor=purple,
    bookmarksnumbered=true,
]{hyperref}

\usepackage[totalwidth=160mm, totalheight=247mm]{geometry}

\usepackage{caption}
\renewcommand{\leq}{\leqslant}
\renewcommand{\geq}{\geqslant}

\newcommand{\bbZ}{\mathbb{Z}}

\newcommand{\bbR}{\mathbb{R}}

\newcommand{\mcA}{\mathcal A}
\newcommand{\mcP}{\mathcal P}
\newcommand{\mcL}{\mathcal L}
\newcommand{\mcU}{\mathcal U}
\newcommand{\mcD}{\mathcal D}

\newcommand{\bfe}{{\bf e}}
\newcommand{\bfv}{{\bf v}}
\newcommand{\bfx}{{\bf x}}
\newcommand{\bfy}{{\bf y}}
\newcommand{\bfP}{{\bf P}}
\newcommand{\bfM}{{\bf M}}
\newcommand{\Pv}{{\mathcal P_{\bf v}}}
\newcommand{\transp}[1]{{{}^\textup{t} #1}}
\newcommand{\EOS}{{{\bf E}_1^*}}
\newcommand{\EOSS}{{{\bf E}_1^*(\sigma)}}

\newcommand{\ub}{{{\bf u}_\beta}}
\newcommand{\vb}{{{\bf v}_\beta}}
\newcommand{\Ms}{{{\bf M}_\sigma}}
\newcommand{\Msinv}{{{\bf M}^{-1}_\sigma}}
\newcommand{\Pc}{{\mathbb P_{\textup c}}}

\newcommand{\dH}{{d_{\textup{H}}}}

\newcommand{\pic}{{\pi_{\textup{c}}}}

\newcommand{\LAR}{\mathcal L_\textup{AR}}

\newcommand{\myvcenter}[1]{\ensuremath{\vcenter{\hbox{#1}}}}

\newtheorem{thm}{Theorem}[section]
\newtheorem{prop}[thm]{Proposition}
\newtheorem{lem}[thm]{Lemma}

\theoremstyle{definition}
\newtheorem{defi}[thm]{Definition}
\newtheorem{exmpl}[thm]{Example}

\title{\textbf{Connectedness of fractals associated with Arnoux-Rauzy substitutions}}

\author[1]{Val\'erie Berth\'e}
\author[1,2]{Timo Jolivet}
\author[3]{Anne Siegel}

\affil[1]{
    LIAFA,
    Universit\'e Paris 7, France
}
\affil[2]{
    FUNDIM,
    Department of Mathematics,
    University of Turku, Finland
}
\affil[3]{
    IRISA,
    Campus de Beaulieu, Rennes, France
}

\date{}

\begin{document}

\maketitle

\begin{abstract}
Rauzy fractals are compact sets with fractal boundary
that can be associated with any unimodular Pisot irreducible substitution.
These fractals can be defined as the Hausdorff limit of a sequence of compact sets,
where each set is a renormalized projection of a finite union of faces of unit cubes.
We exploit this combinatorial definition to prove the connectedness
of the Rauzy fractal associated with any finite product of three-letter Arnoux-Rauzy substitutions.
\end{abstract}

\section{Introduction}
\begin{figure}[ht]
\centering
\begin{tabular}{ccc}
    \myvcenter{\includegraphics[width=30mm]{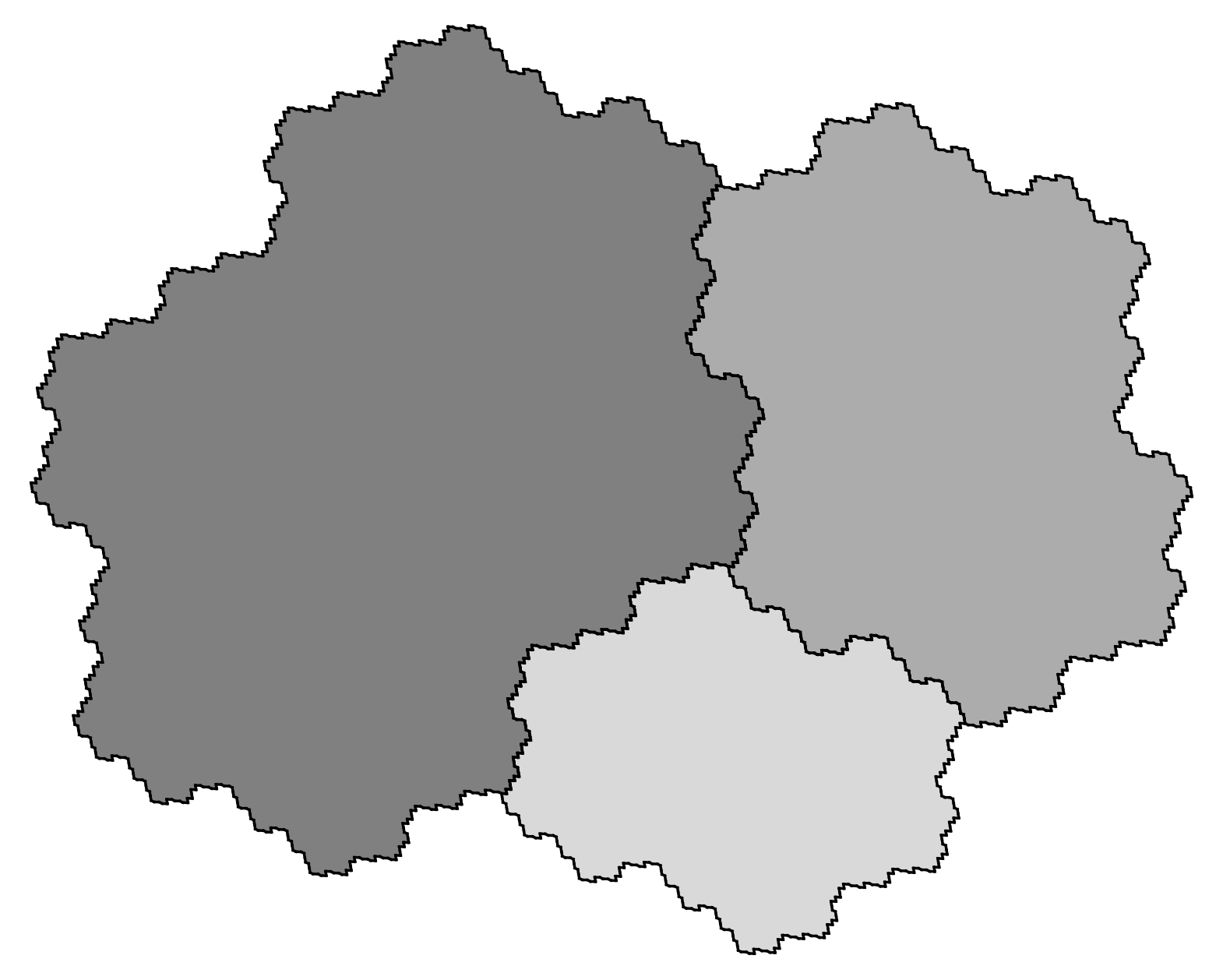}} &
    \hspace{5mm} \myvcenter{$\longmapsto$} \hspace{5mm} &
    \myvcenter{\includegraphics[width=30mm]{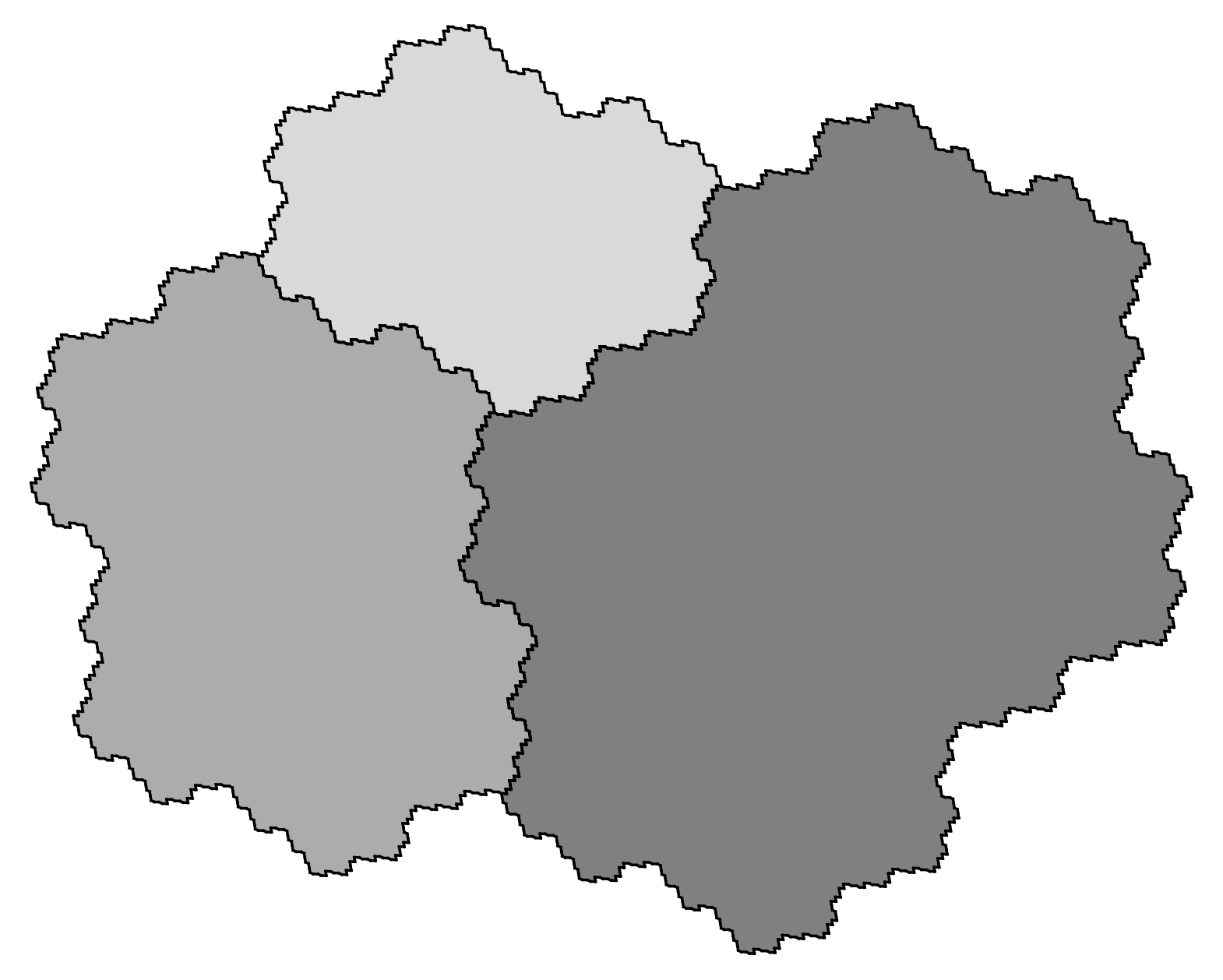}}
\end{tabular}
\caption{Domain exchange in the Tribonacci fractal.}
\label{fig:triboexchange}
\end{figure}

Rauzy fractals are compact sets with fractal boundary
that can be associated with any unimodular Pisot irreducible substitution.
(See Definition~\ref{def:Pisot} and~\ref{def:RF} for precise definitions.)
They first appeared in the work of Rauzy~\cite{Rau82},
who generalized the theory of interval exchange transformations
by defining a \emph{domain} exchange transformation
of three pieces in $\bbR^2$.
Each of the three pieces is translated along a vector
(one distinct vector for each piece)
in order to give a different partition of the same shape
(see Figure~\ref{fig:triboexchange}).
These fractals were also discussed in the later work of Thurston~\cite{Thu89}
in the context of numeration systems in non-integer bases.

\paragraph{Properties of Rauzy fractals}
Originally, Rauzy considered the
\emph{Tribonacci substitution} $\sigma: 1 \mapsto 12, 2 \mapsto 13, 3 \mapsto 1$,
and he proved in \cite{Rau82} that the dynamics of the symbolic dynamical system generated by $\sigma$
is realized by the domain exchange in the Rauzy fractal associated with $\sigma$.
Rauzy's results have then been generalized by Arnoux and Ito~\cite{AI01} as follows:
for any unimodular Pisot irreducible substitution $\sigma$ on $d$ letters,
the subshift $X_\sigma$ generated by $\sigma$ is semi-conjugate
to a domain exchange of $d$ pieces in the Rauzy fractal associated with $\sigma$ (a compact subset of $\bbR^{d-1}$),
under the assumption that $\sigma$ verifies a combinatorial condition called
the \emph{strong coincidence condition} (see also~\cite{CS01}).

Rauzy fractals also provide an explicit way to prove that
the subshift $X_\sigma$ is semi-conjugate to a translation on the $(d-1)$-dimensional torus
if $\sigma$ satisfies the \emph{super coincidence condition} introduced in~\cite{IR06, BK06}.
It is currently not known if a Pisot irreducible substitution $\sigma$ always verifies the
strong and the super coincidence conditions;
the Pisot conjecture states that this is always the case (see \cite{Que10,PF02,CANT}).

Let us mention some other properties of Rauzy fractals.
In numeration systems, they provide natural extensions
of $\beta$-transformations with relevant algebraic properties \cite{ABBS08}.
In theoretical physics, they are good candidates for explicit
cut-and-project schemes which model quasicrystals~\cite{GVG04}.
In discrete geometry, they are related to
discrete plane generation via
multidimensional continued fraction algorithms~\cite{IO93, ABFJ07, Fer09}.

\begin{figure}[ht]%
\centering
\subfloat[][$1 \mapsto 12, 2 \mapsto 31, 3 \mapsto 1$]{%
    \label{fig:topocool-a}%
    \includegraphics[height=32mm]{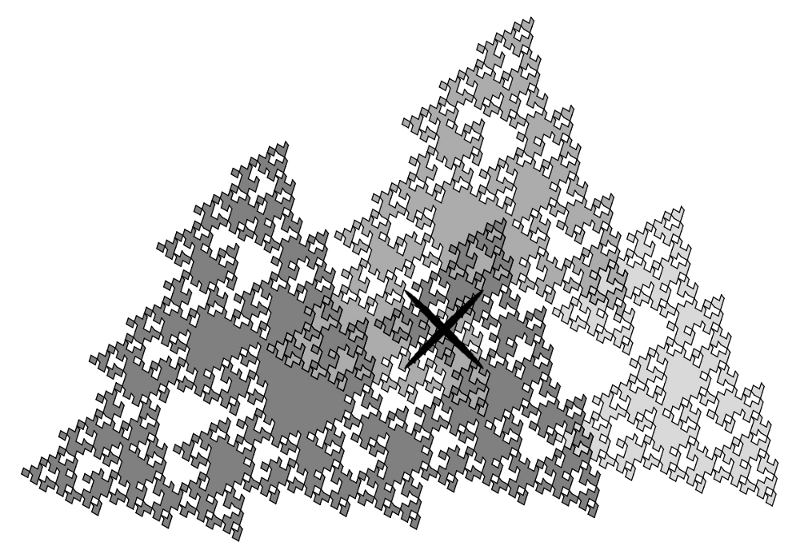}
} \hfil
\subfloat[][$1 \mapsto 3, 2 \mapsto 23, 3 \mapsto 31223$]{%
    \label{fig:topocool-b}%
    \includegraphics[height=32mm]{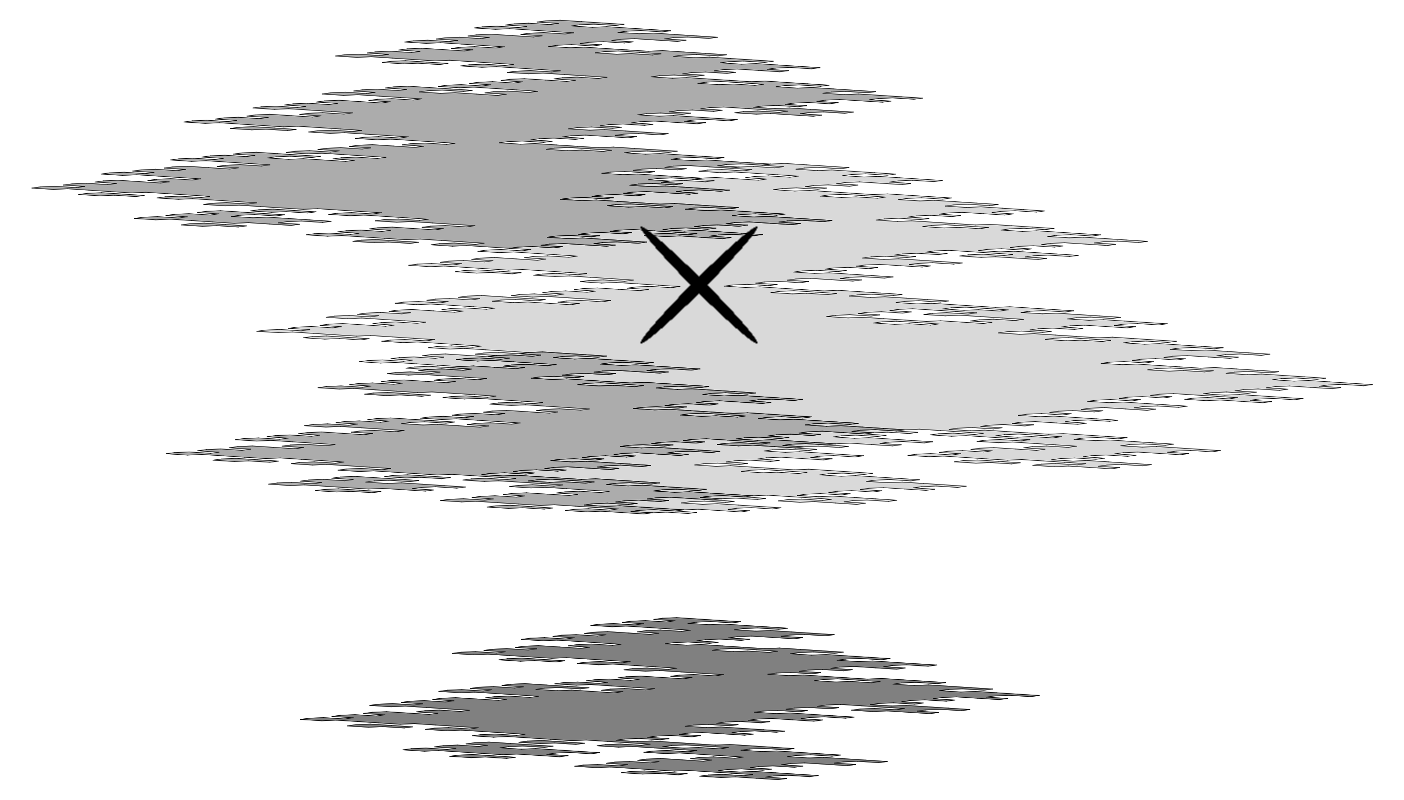}
}\\
\subfloat[][$1 \mapsto 21111, 2 \mapsto 31111, 3 \mapsto 1$]{%
    \label{fig:topocool-c}%
    \includegraphics[height=25mm]{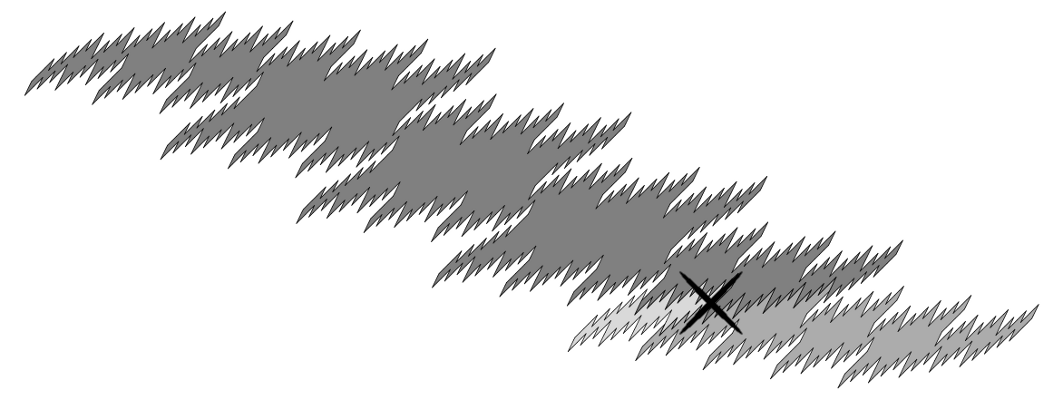}
} \hfil
\subfloat[][$1 \mapsto 123, 2 \mapsto 1, 3 \mapsto 31$]{%
    \label{fig:topocool-d}%
    \includegraphics[height=30mm]{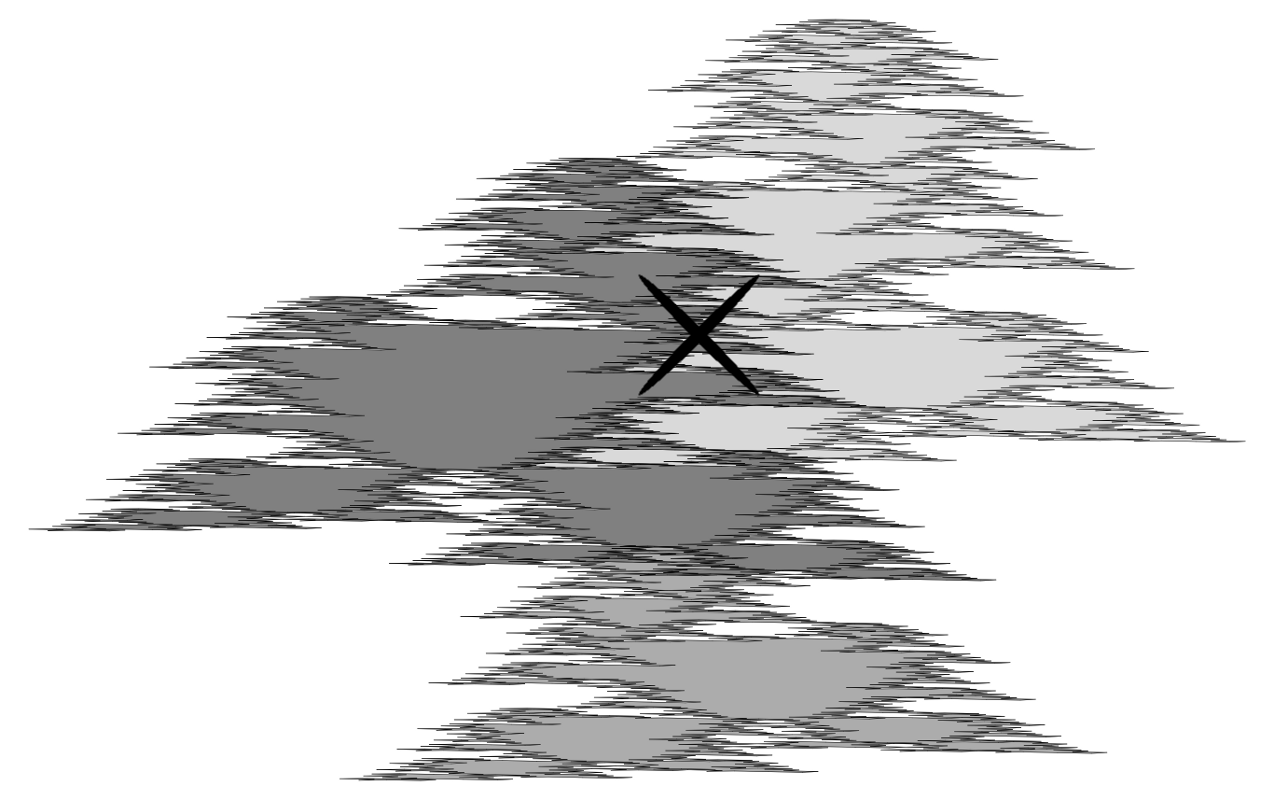}
}%
\caption[]{Examples of Rauzy fractals:
\subref{fig:topocool-a} is connected and has uncountable fundamental group;
\subref{fig:topocool-b} is disconnected and $\mathbf 0$ is not an inner point;
\subref{fig:topocool-c} is connected, has uncountable fundamental group, and $\mathbf 0$ is an inner point;
\subref{fig:topocool-d} is connected and $\mathbf 0$ is not an inner point.
Black crosses mark the origin.
}%
\label{fig:topocool}%
\end{figure}

\paragraph{Topology of Rauzy fractals}
The Rauzy fractal associated with the Tribonacci substitution
(Figure~\ref{fig:triboexchange}) enjoys good topological properties:
the origin is an inner point, and it is homeomorphic to a closed disc~\cite{Mes00}.
However, the topology of Rauzy fractals can be very complicated in general:
they can fail to be connected or simply connected,
and the origin is not always an inner point of the set as shown in Figure~\ref{fig:topocool}.
More examples are given in~\cite{ST10}.

The study of topological properties of Rauzy fractals has many applications.
The connectedness of the associated Rauzy fractal implies
the connectedness of the domains used to provide geometrical interpretations of the dynamics of substitutions:
the fundamental domains used to encode the domain exchange transformation or the toral translation,
and the Markov partition of the toral automorphism provided by
the incidence matrix of the underlying substitution.
This is stated in details in Theorem~\ref{thm:ARmarkov}.

Moreover, in Diophantine approximation, information on the size of the largest ball contained in the Rauzy fractal
allows the determination of the sequence of best approximations, with respect to a specific norm, for some two-dimensional vectors
provided by non-totally real cubic Pisot units \cite{HM06}.

In number theory,
finite greedy expansions in non-integer bases ($\beta$-numeration systems)
are closely related to the inner points of the fractal,
and the connectedness of the fractal is conjectured to guarantee explicit relations between
the norm of $\beta$ and the $\beta$-expansion of $1$~\cite{AN05}.
The properties of rational numbers with purely periodic $\beta$-expansions are
closely related to the shape of the boundary of the Rauzy fractal~\cite{ABBS08, AFSS10}.
In discrete geometry, studying the position of the origin in the fractal
allows us to study the structure of discrete planes and to generate them~\cite{IO93, BLPP13}.

Lastly, cut-points of the fractal are related to some topological invariants of tiling spaces \cite{BDS09}.

\paragraph{Arnoux-Rauzy substitutions}
Sturmian sequences are a classical object of symbolic dynamics.
They are the infinite sequences of two letters with factor complexity $n+1$
(\emph{i.e.}, they have exactly $n+1$ factors of length $n$),
and they correspond to natural codings of irrational rotations on the circle~\cite{MH40}.
These sequences are also closely related to continued fraction expansions of real numbers;
see~\cite[Chap. 6]{PF02} and~\cite[Chap. 2]{Lot97}
for a detailed survey of their properties.

Arnoux-Rauzy sequences were introduced in~\cite{AR91} to generalize
Sturmian sequences to $3$-letter alphabets.
They are the infinite sequences of three symbols obtained by
iterating the \emph{Arnoux-Rauzy substitutions}
$\sigma_1$, $\sigma_2$, $\sigma_3$ (defined in Section \ref{sec:ARsub})
where each $\sigma_i$ occurs infinitely often.
These sequences have factor complexity $2n+1$ and they yield an algorithm for simultaneous approximations of
some particular pairs of algebraic numbers \cite{CFM08}.
It was conjectured that Arnoux-Rauzy sequences correspond to natural codings
of translations on the two-dimensional torus (as in the Sturmian case with rotations on the circle),
but this conjecture has been disproved in~\cite{CFZ00}.
For more references about Arnoux-Rauzy sequences, see~\cite{BFZ05, CC06}.

\paragraph{Our results}
In the case where $\sigma$ is a $2$-letter substitution,
topological properties of Rauzy fractals (which are subsets of $\bbR$) are fully understood:
the fractal is connected (\emph{i.e.}, an interval) if and only if the substitution is Sturmian~\cite{EI98},~\cite[Chap. 9]{PF02}.
Several conditions equivalent to the latter property in the $2$-letter case are summarized in~\cite{BFS12}.

In this article, we study the infinite family of the Rauzy fractals associated with
finite products of $3$-letter Arnoux-Rauzy substitutions,
and we prove that they are connected (Theorem~\ref{thm:ARfractal}).

To this end, we use a combinatorial characterization
of Rauzy fractals given by Arnoux and Ito in~\cite{AI01} (see Definition~\ref{def:RF} below).
For any unimodular Pisot irreducible substitution $\sigma$ on $d$ letters,
they define a \emph{dual substitution} $\EOSS$,
that does not act words, but on faces of unit cubes in $\bbR^d$.
The Rauzy fractal associated with $\sigma$
can then be obtained by iterating $\EOSS$
starting from a small set of unit faces
(for example \myvcenter{\begin{tikzpicture}
[scale=0.8, x={(-0.216506cm,-0.125000cm)}, y={(0.216506cm,-0.125000cm)}, z={(0.000000cm,0.250000cm)}]
\definecolor{facecolor}{rgb}{0.800000,0.800000,0.800000}
\fill[fill=facecolor, draw=black, shift={(0,0,0)}]
(0, 0, 0) -- (0, 1, 0) -- (0, 1, 1) -- (0, 0, 1) -- cycle;
\fill[fill=facecolor, draw=black, shift={(0,0,0)}]
(0, 0, 0) -- (0, 0, 1) -- (1, 0, 1) -- (1, 0, 0) -- cycle;
\fill[fill=facecolor, draw=black, shift={(0,0,0)}]
(0, 0, 0) -- (1, 0, 0) -- (1, 1, 0) -- (0, 1, 0) -- cycle;
\end{tikzpicture}}
 when $d=3$).
This gives an increasing sequence of finite sets of unit faces in $\bbR^d$ which,
if projected on a particular hyperplane and renormalized appropriately at each step,
admits a Hausdorff limit that is equal to the Rauzy fractal.

Our results are based on an alternative description
of $\EOSS$ substitutions, introduced in~\cite{IO93, IO94, ABI02}.
It consists in trying to compute the image $\EOSS(P)$ of a set of unit faces $P$
by \emph{concatenating} the images of the elements of $P$,
instead of using the definition of $\EOSS$
to compute the new position of each face.
This is similar to the concatenation relation $\sigma(uv) = \sigma(u)\sigma(v)$
which is valid for words but difficult to generalize to higher dimensions.

A proof of Theorem~\ref{thm:ARfractal} has been announced in~\cite{Can03},
relying on completely different methods, but it has not been published.

Let us note that the results published in \cite{BJS12}
use techniques similar to that of the current article,
but with a different goal,
namely that the $\EOSS^n(\myvcenter{\begin{tikzpicture}
[scale=0.8, x={(-0.216506cm,-0.125000cm)}, y={(0.216506cm,-0.125000cm)}, z={(0.000000cm,0.250000cm)}]
\definecolor{facecolor}{rgb}{0.800000,0.800000,0.800000}
\fill[fill=facecolor, draw=black, shift={(0,0,0)}]
(0, 0, 0) -- (0, 1, 0) -- (0, 1, 1) -- (0, 0, 1) -- cycle;
\fill[fill=facecolor, draw=black, shift={(0,0,0)}]
(0, 0, 0) -- (0, 0, 1) -- (1, 0, 1) -- (1, 0, 0) -- cycle;
\fill[fill=facecolor, draw=black, shift={(0,0,0)}]
(0, 0, 0) -- (1, 0, 0) -- (1, 1, 0) -- (0, 1, 0) -- cycle;
\end{tikzpicture}}
)$ cover arbitrarily large discs.
As a consequence, it is proved in \cite{BJS12} that the symbolic dynamical system generated by
a purely substitutive Arnoux-Rauzy sequence is measurably conjugate to a toral translation.

\paragraph{Outline of the paper}
In Section~\ref{sec:prelim}, we give the definition of
dual substitutions and relate them with discrete planes,
we give a definition of the Rauzy fractal,
and we introduce Arnoux-Rauzy substitutions.
In Section~\ref{sec:cover}, we establish a combinatorial sufficient
condition for the connectedness of Rauzy fractals.
In Section~\ref{sec:appl}, we apply the results of Section~\ref{sec:cover}
to prove the connectedness of the fractal associated with an
Arnoux-Rauzy substitution (Theorem~\ref{thm:ARapprox} and~\ref{thm:ARfractal}).
We deduce the existence of a connected Markov partition
for the associated toral automorphisms (Theorem~\ref{thm:ARmarkov}).
Finally, in Section~\ref{sec:counterex}, we provide examples to show that
some possible generalizations of Theorem~\ref{thm:ARfractal} do not hold.

\section{Preliminaries}
\label{sec:prelim}
The definitions and results of this section are valid in any dimension
but we state them for dimension $3$ only, since our main results
(Theorem~\ref{thm:ARapprox} and~\ref{thm:ARfractal})
hold in dimension $3$.

\subsection{Discrete planes and unit faces}
We start by giving a geometric definition of discrete planes,
following~\cite{Rev91, IO93, IO94, ABI02}.
Let $(\bfe_1, \bfe_2, \bfe_3)$ denote the canonical basis of $\bbR^3$
and let $\bbR_{>0}^3$ denote the set of vectors of $\bbR^3$ with positive entries.
Recall that the plane of (non-zero) normal vector $\bfv \in \bbR^3$
is the set of points $\bfx \in \bbR^3$ such that $\langle \bfx, \bfv \rangle = 0$.

\begin{defi}[Discrete plane]
\label{def:discreteplane}
Let $\bfv \in \bbR_{>0}^3$ and
let $\mathcal S$ be the union of the unit cubes with integer coordinates
that intersect the lower half-space $\{\bfx \in \bbR^3 : \langle \bfx, \bfv \rangle < 0\}$.
The \emph{discrete plane $\Pv$ of normal vector $\bfv$} is the boundary of $\mathcal S$.
\end{defi}

\begin{figure}[ht]
\centering
\includegraphics[width=0.9\linewidth]{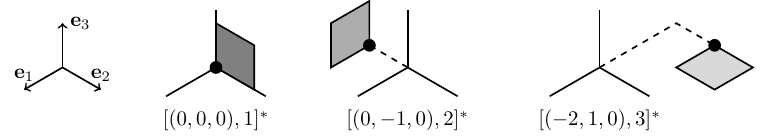}
\caption{Three unit faces of different types.}
\label{fig:faces}
\end{figure}

A discrete plane can be seen as a union of unit faces of three different types.
Let $i \in \{1,2,3\}$ and $\bfx \in \bbZ^3$.
The \emph{unit face type $i$ at point $\bfx$} is the set $[\bfx, i]^\star$ defined by
\[
\begin{array}{ccc}
  \,[\bfx, 1]^\star & = & \{\bfx + \lambda \bfe_2 + \mu \bfe_3 : \lambda, \mu \in [0,1]\} \\
  \,[\bfx, 2]^\star & = & \{\bfx + \lambda \bfe_1 + \mu \bfe_3 : \lambda, \mu \in [0,1]\} \\
  \,[\bfx, 3]^\star & = & \{\bfx + \lambda \bfe_1 + \mu \bfe_2 : \lambda, \mu \in [0,1]\}
\end{array}
\]
(see Figure~\ref{fig:faces}).
The type $i$ of face $[\bfx, i]^\star$ corresponds to the canonical vector $\bfe_i$
to which it is orthogonal.
If $\bfx \in \bbZ^3$, we can write $\bfy + [\bfx, i]^\star$ instead of $[\bfx + \bfy, i]^\star$,
and we denote by $\bfx + X$ the translation of a union of faces $X$ by $\bfx$.
Let us remark that Definition~\ref{def:discreteplane} implies that the set
$[\mathbf 0, 1]^\star \cup [\mathbf 0, 2]^\star \cup [\mathbf 0, 3]^\star$ is included in every discrete plane.

The following proposition gives an alternative definition of discrete planes,
where the belonging of each face to a plane is translated into
an inequality on scalar products.

\begin{prop}[\cite{ABI02,ABS04}]
Let $\bfv \in \bbR_{>0}^3$.
The discrete plane $\Pv$ is the union of faces $[\bfx, i]^\star$
satisfying $0 \leq \langle \bfx, \bfv \rangle < \langle \bfe_i, \bfv \rangle$.
\end{prop}

\subsection{Dual substitutions}
We start with the classical notion of unidimensional substitution.
Let $\mathcal A = \{1, \ldots, d\}$ be a finite alphabet,
and denote by $\mcA^\star$ the set of finite words over $\mcA$.
A~\emph{substitution} is a function
$\sigma : \mathcal A^\star \rightarrow \mathcal A^\star$
such that $\sigma(uv) = \sigma(u)\sigma(v)$ for all words $u, v \in \mathcal A^\star$,
and such that the image of each letter of $\mathcal A$ is not the empty word.
The \emph{incidence matrix} $\Ms$ of $\sigma$
is the square matrix of size $d \times d$
defined by $\Ms = (m_{ij})$,
where $m_{ij}$ is the number of occurrences of the letter $i$ in $\sigma(j)$.
A substitution is \emph{unimodular} if the $\det(\bfM_\sigma) = \pm 1$.

A classical example of a substitution is the Tribonacci substitution
introduced by Rauzy in~\cite{Rau82},
whose action on $\{1,2,3\}^\star$ and incidence matrix are given by
\[
\sigma \ : \ \left\{
\begin{array}{rcl}
1 & \mapsto & 12 \\
2 & \mapsto & 13 \\
3 & \mapsto & 1
\end{array}
\right.
\qquad \text{ and } \qquad
\Ms \ = \
\begin{pmatrix}
  1 & 1 & 1 \\
  1 & 0 & 0 \\
  0 & 1 & 0
\end{pmatrix}.
\]

We now introduce \emph{dual substitutions},
which act not on unidimensional words, but on unit faces in $\bbR^3$.
This formalism was sketched by Ito and Ohtsuki~\cite{IO93, IO94},
and then refined later by Arnoux and Ito in~\cite{AI01} (see Definition~\ref{def:gensub} below),
where they also highlight the connections between dual substitutions
and discrete planes (Propositions~\ref{prop:distinctfaces} and~\ref{prop:imgplane} below).

\begin{defi}[Dual substitution]
\label{def:gensub}
Let $\sigma : \{1,2,3\}^\star \rightarrow \{1,2,3\}^\star$ be a unimodular substitution.
The \emph{dual substitution} associated with $\sigma$,
denoted by $\EOSS$, is defined by
\[
\EOSS([\bfx, i]^\star) \ = \
\bigcup_{j = 1,2,3}
\;
\bigcup_{s | \sigma(j) = pis}
[\Msinv(\bfx + \bfP(s)), j]^\star,
\]
where $\bfP : \mcA^\star \rightarrow \bbZ^n$
is the \emph{Abelianization map} defined by
$\bfP(w) = (|w|_1, |w|_2, |w|_3)$,
where $|w|_i$ denotes the number of occurrences of $i$ in $w$.
We extend this definition to any union of unit faces:
\[
\EOSS(X_1 \cup X_2)
  \ = \ \EOSS(X_1) \cup \EOSS(X_2).
\]
\end{defi}

Let us remark that $\EOSS$ is completely described by $\Ms$ and
the images of the faces $[\mathbf 0, 1]^\star$, $[\mathbf 0, 2]^\star$ and $[\mathbf 0, 3]^\star$,
because $\EOSS([\bfx, i]^\star) = \Msinv \bfx + \EOS([\mathbf 0, i]^\star)$
for every unimodular substitution $\sigma$ and every face $[\bfx, i]^\star$.
It is also worth noticing that
$\EOS(\sigma \circ \sigma') = \EOS(\sigma') \circ \EOS(\sigma)$
holds for all unimodular substitutions $\sigma$ and $\sigma'$;
see~\cite{AI01} for more details.

\begin{exmpl}
Let $\sigma$ be the Tribonacci substitution
$1 \mapsto 12$, $2 \mapsto 13$, $3 \mapsto 1$.
The action of $\EOSS$ on unit faces is given by
\[
\begin{array}{rcl}
\EOSS([\bfx, 1]^\star) & = &
  \Msinv \bfx + [(1,0,-1), 1]^\star \cup [(0,1,-1),2]^\star \cup [\mathbf 0, 3]^\star\\
\EOSS([\bfx, 2]^\star) & = &
  \Msinv\bfx + [\mathbf 0, 1]^\star\\
\EOSS([\bfx, 3]^\star) & = &
  \Msinv\bfx + [\mathbf 0, 2]^\star,
\end{array}
\]
which can be represented graphically as follows:
\[
\includegraphics[width=0.9\linewidth]{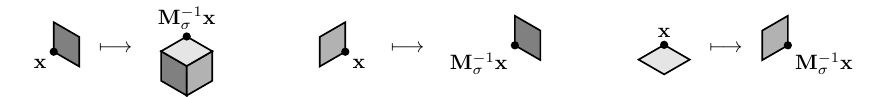}
\]
\end{exmpl}

In general, the images of two distinct faces are not necessarily disjoint,
but it is the case when the faces belong to a common discrete plane,
as stated in Proposition~\ref{prop:distinctfaces}.

\begin{prop}[\cite{AI01}]
\label{prop:distinctfaces}
Let $\bfv \in \bbR_{>0}^3$.
If $[\bfx, i]^\star$ and $[\bfx', i']^\star$ are two distinct faces of $\Pv$,
then the sets $\EOSS([\bfx, i]^\star)$ and $\EOSS([\bfx', i']^\star)$ are disjoint up to a set of measure zero.
\end{prop}

Proposition~\ref{prop:imgplane} states that the image of a discrete plane by
a dual substitution is again a discrete plane.

\begin{prop}[\cite{AI01, Fer06}]
\label{prop:imgplane}
Let $\bfv \in \bbR_{>0}^3$
and $\sigma$ be a unimodular substitution.
We have $E_1^\star(\sigma)(\Pv) = \mcP_{\transp \Ms \bfv}$.
\end{prop}

\subsection{The Rauzy fractal associated with a substitution}
We will now give the definition of the Rauzy fractal associated with
a unimodular Pisot irreducible substitution, as in~\cite{AI01}.
We recall that a \emph{Pisot number} is a real algebraic integer
greater than $1$ whose conjugates have absolute value less than $1$.

\begin{defi}[Pisot irreducible substitution]
\label{def:Pisot}
A unimodular substitution $\sigma$ is \emph{Pisot irreducible}
if the maximal eigenvalue of $\Ms$ is a Pisot number,
and if the characteristic polynomial of $\Ms$ is irreducible.
\end{defi}

Let $\sigma : \{1,2,3\}^\star \rightarrow \{1,2,3\}^\star$ be a unimodular Pisot irreducible substitution,
and let $\beta$ be the maximal eigenvalue of $\Ms$.
We denote by $\ub$ a left-eigenvector of $\Ms$ associated with $\beta$,
and by $\vb$ a right-eigenvector of $\Ms$ associated with $\beta$ (with positive coordinates).
Such a vector $\vb$ indeed exists, thanks to the Perron-Frobenius theorem
and the fact that the matrix of a unimodular Pisot irreducible substitution
is always primitive; see~\cite{CS01}.

\begin{defi}[Contracting plane]
Let $\sigma$ be a unimodular Pisot irreducible substitution.
The \emph{contracting plane} $\Pc$ associated with $\sigma$
is the plane of normal vector $\vb$.
We denote by $\pic : \bbR^3 \rightarrow \Pc$ the projection
from $\bbR^3$ to $\Pc$ along $\ub$.
\end{defi}

Let $\mcU = [{\bf 0},1]^\star \cup [{\bf 0},2]^\star \cup [{\bf 0},3]^\star$.
Proposition~\ref{prop:imgplane} enables us to iterate
$\EOSS$ on $\mcU$ in order to obtain an infinite sequence of patterns of increasing size
that are included in the discrete plane of normal vector $\vb$,
\emph{i.e.}, the discretization of the contracting plane.
Let us remark that the set $\mcP_\vb$ is indeed a discrete plane, because $\vb$ has positive coordinates.

It is possible to project and renormalize the patterns by applying $\Ms \circ \pic$,
in order to obtain a sequence of subsets of the contracting plane $\Pc$
that converges to a compact subset of $\Pc$.
More precisely, for $n \geq 0$, let
\[
\mcD_n \ = \ {\bf M}^{n}_\sigma \circ \pic  \circ \EOSS^n(\mcU).
\]
Arnoux and Ito proved in~\cite{AI01} that $(\mcD_n)_{n \geq 0}$ is a convergent sequence
in the metric space of compact subsets of $\Pc$ together with the Hausdorff distance,
as illustrated below.
\[
    \includegraphics[scale=0.08]{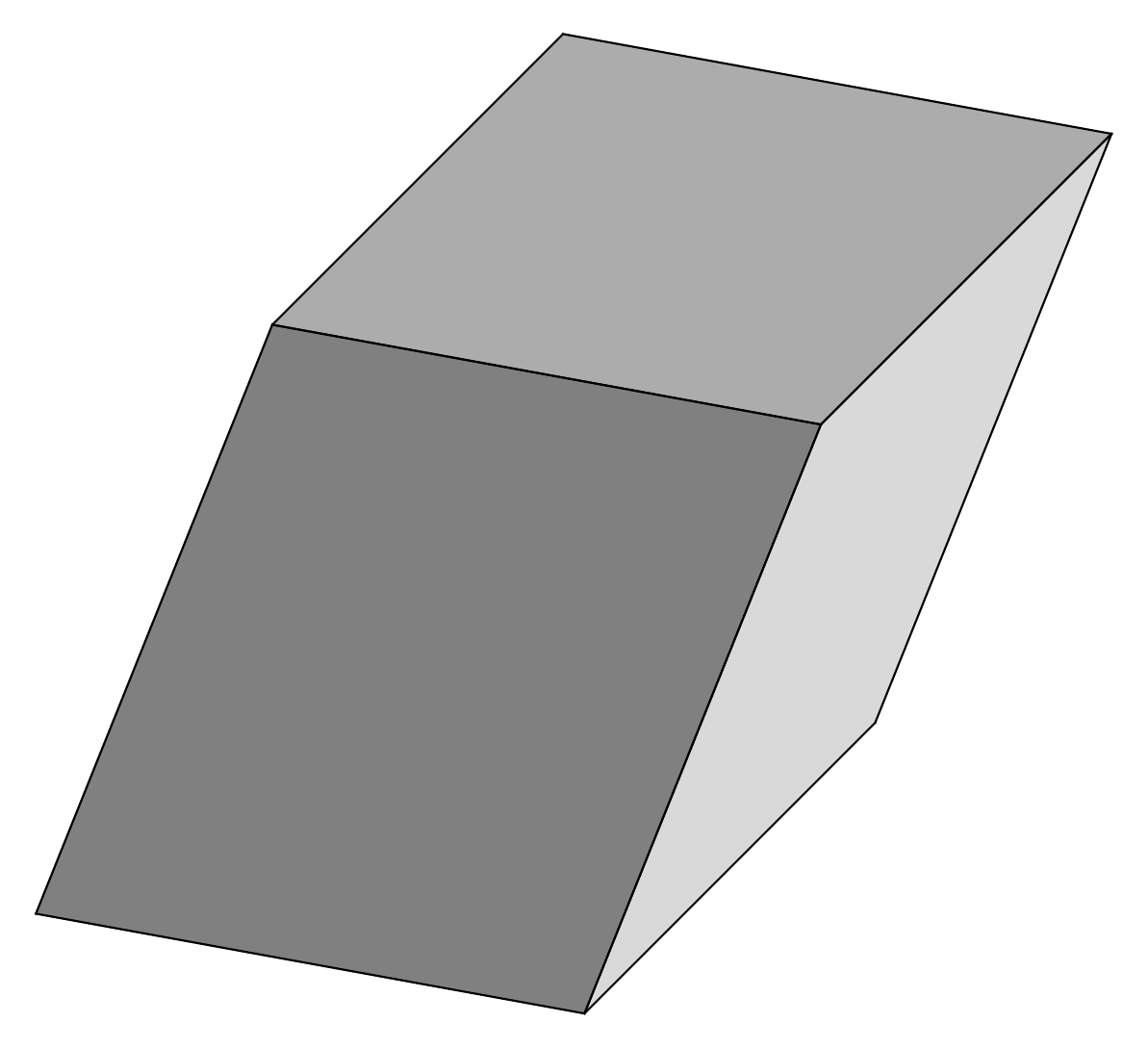} \quad
    \includegraphics[scale=0.08]{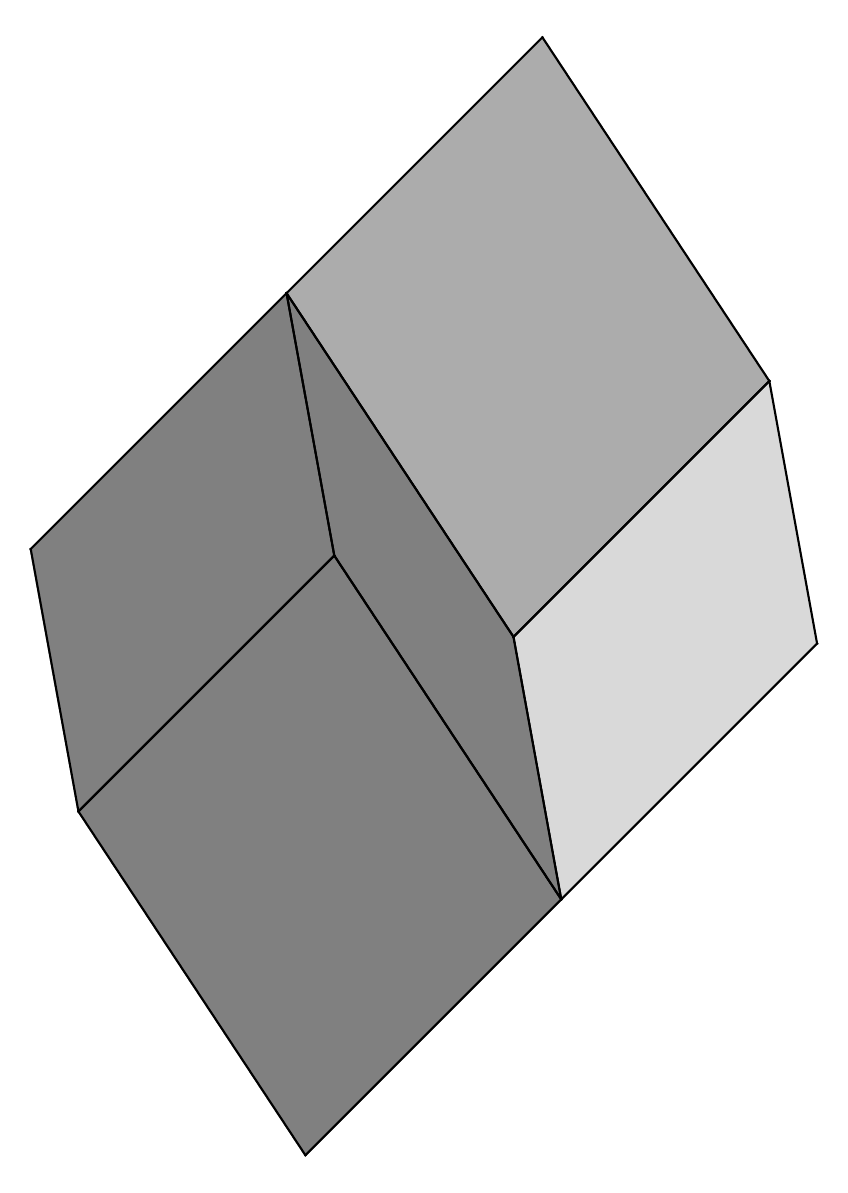} \quad
    \includegraphics[scale=0.08]{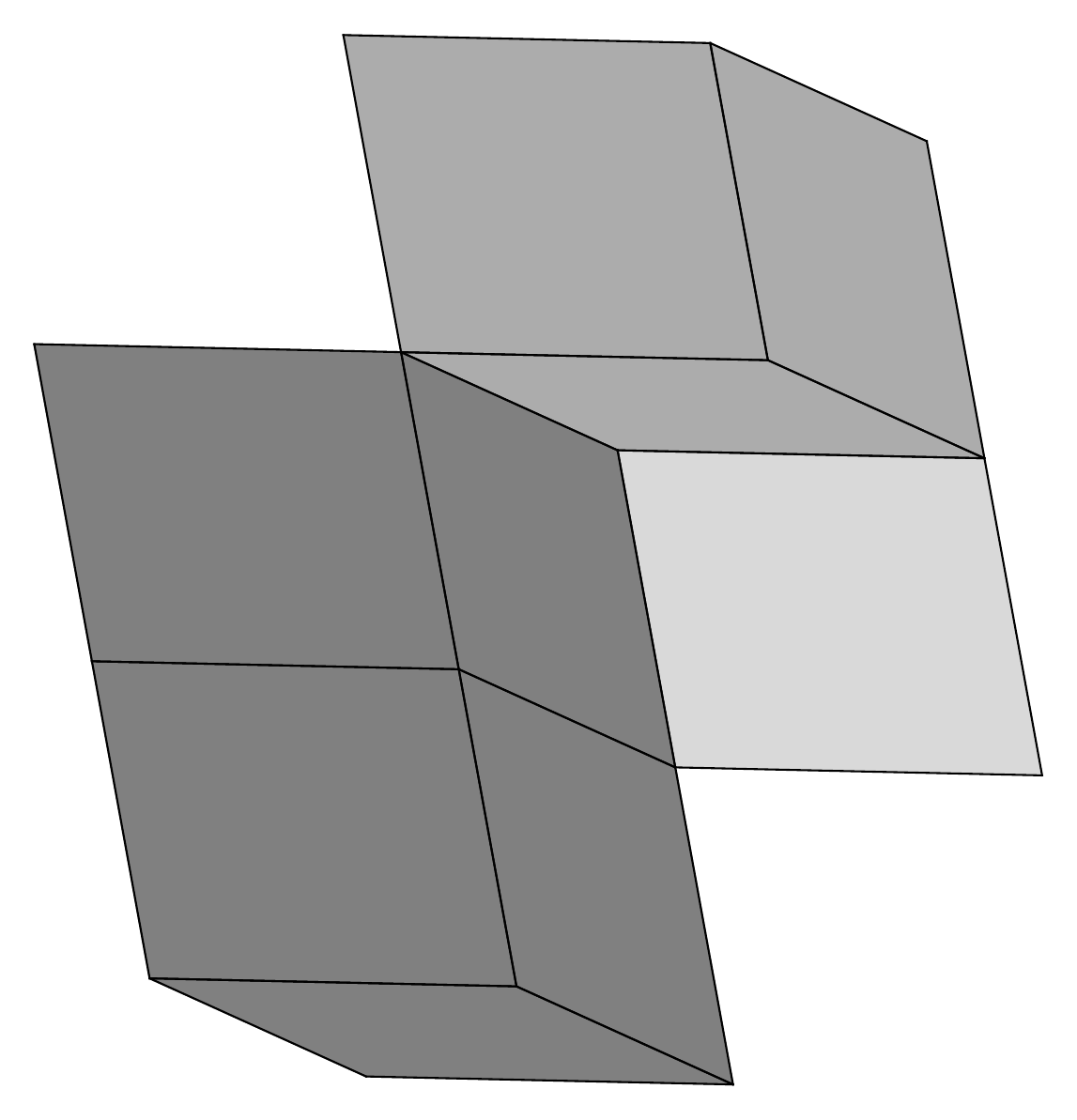} \quad
    \includegraphics[scale=0.08]{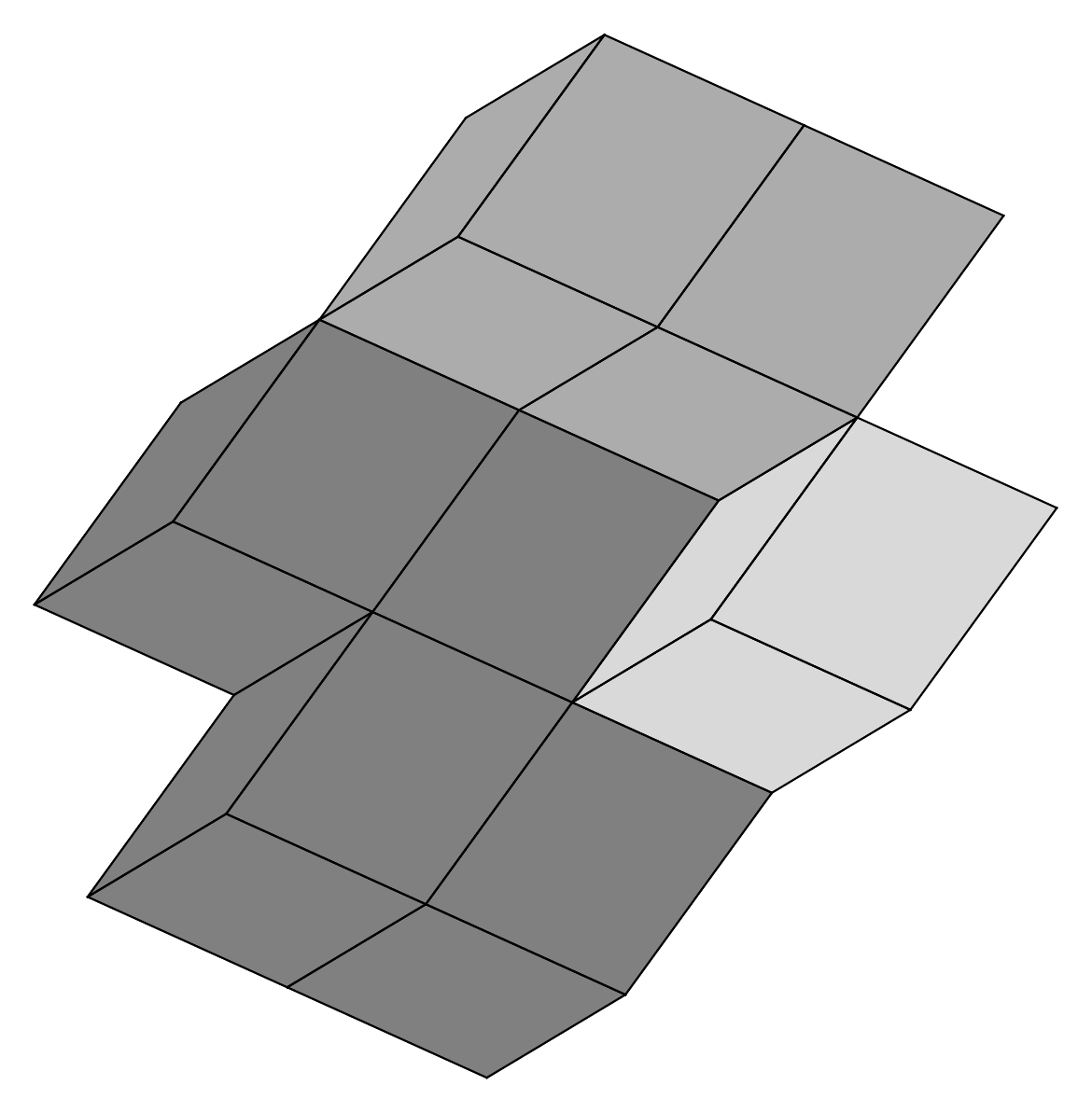} \quad
    \includegraphics[scale=0.08]{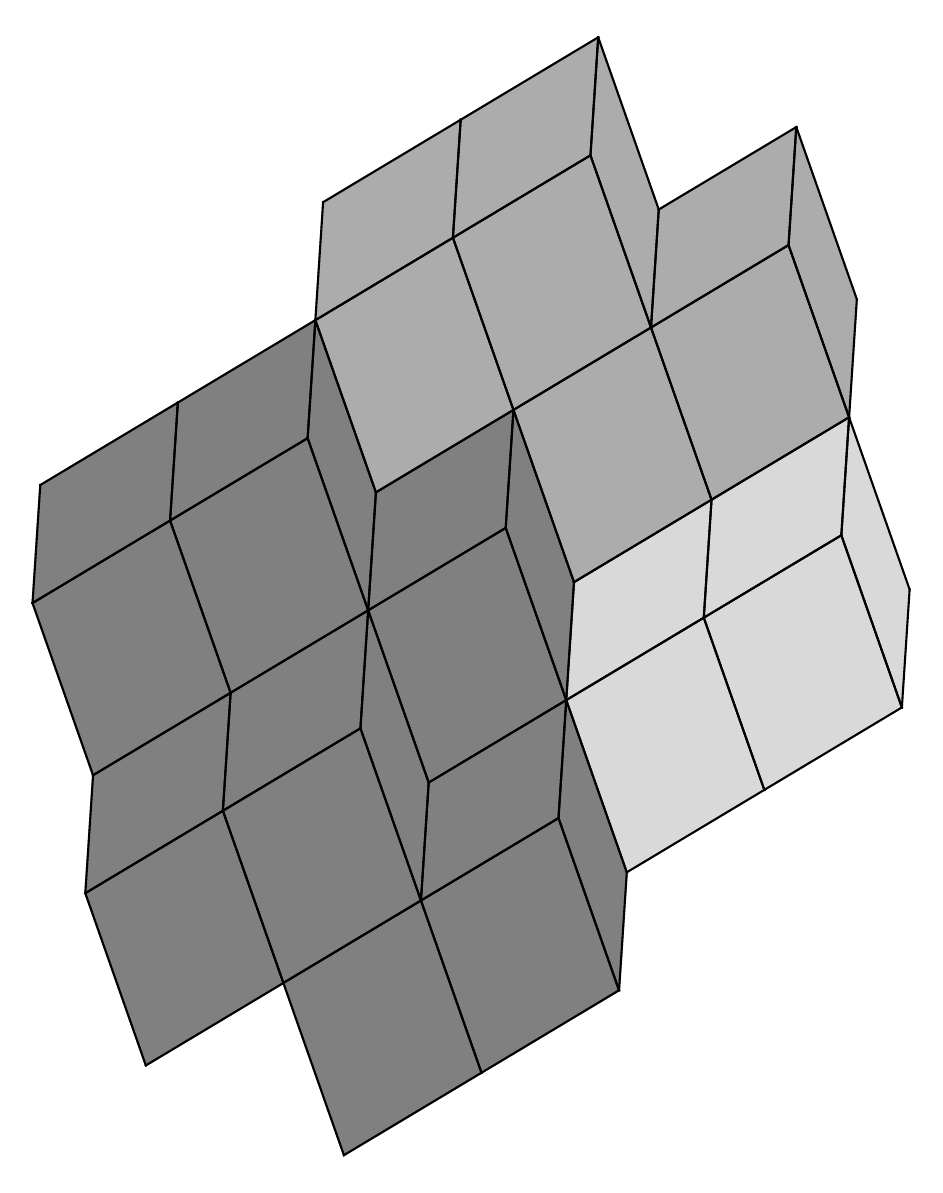} \quad
    \includegraphics[scale=0.08]{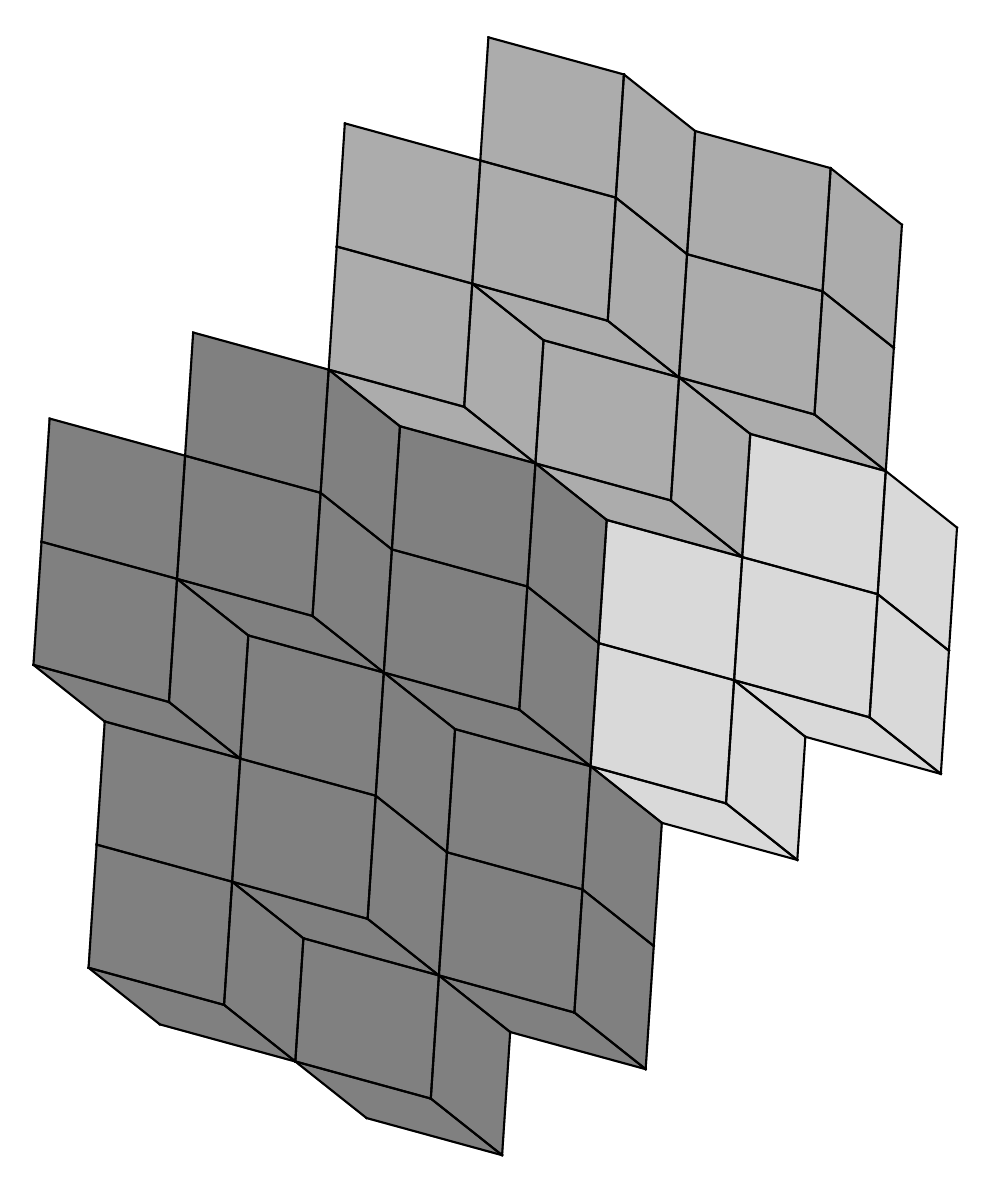} \quad
    \includegraphics[scale=0.08]{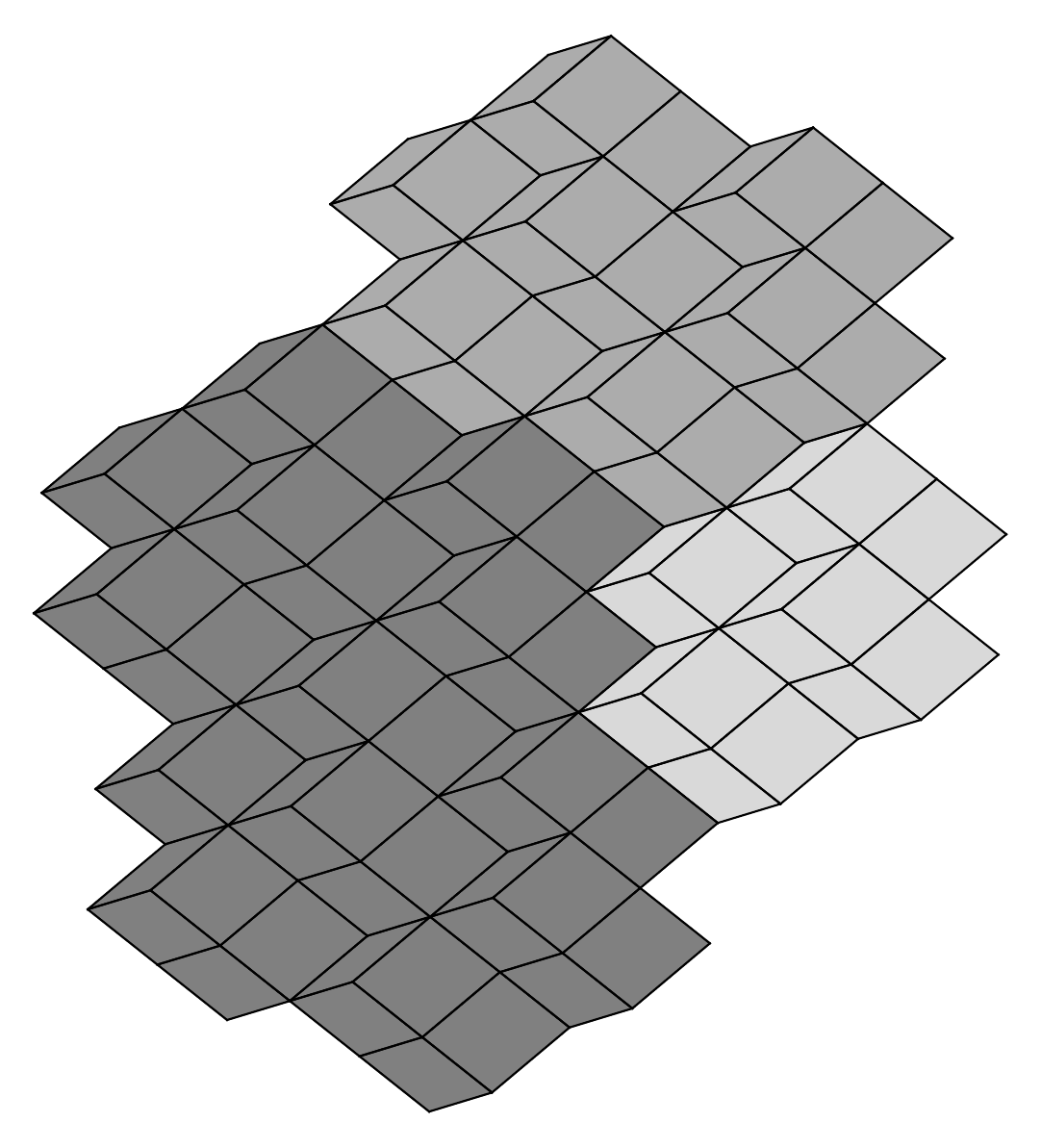}
\]

\begin{defi}[Rauzy fractal]
\label{def:RF}
Let $\sigma : \{1,2,3\}^\star \rightarrow \{1,2,3\}^\star$ be a unimodular Pisot irreducible substitution.
The \emph{Rauzy fractal} associated with $\sigma$ is the Hausdorff limit
of the sequence $(\mcD_n)_{n \geq 0}$.
\end{defi}

\subsection{Arnoux-Rauzy substitutions}
\label{sec:ARsub}
Let $\sigma_1$, $\sigma_2$, $\sigma_3$ be the \emph{Arnoux-Rauzy substitutions}~\cite{AR91} defined by
\[
\sigma_1 \ : \ \left\{
  \begin{array}{rcl}
  1 &\mapsto& 1 \\
  2 &\mapsto& 21 \\
  3 &\mapsto& 31
  \end{array}
\right.,
\qquad
\sigma_2 \ : \ \left\{
  \begin{array}{rcl}
  1 &\mapsto& 12 \\
  2 &\mapsto& 2 \\
  3 &\mapsto& 32
  \end{array}
\right.,
\qquad
\sigma_3 \ : \ \left\{
  \begin{array}{rcl}
  1 &\mapsto& 13 \\
  2 &\mapsto& 23 \\
  3 &\mapsto& 3
  \end{array}
\right..
\]

The following proposition enables us to define Rauzy fractals associated with
finite products of Arnoux-Rauzy substitutions.
\begin{prop}[\cite{AI01}]
A finite product of Arnoux-Rauzy substitutions where each $\sigma_i$ appears at least once
is a unimodular Pisot irreducible substitution.
\end{prop}

\begin{figure}[ht]%
\centering
\includegraphics[height=35mm]{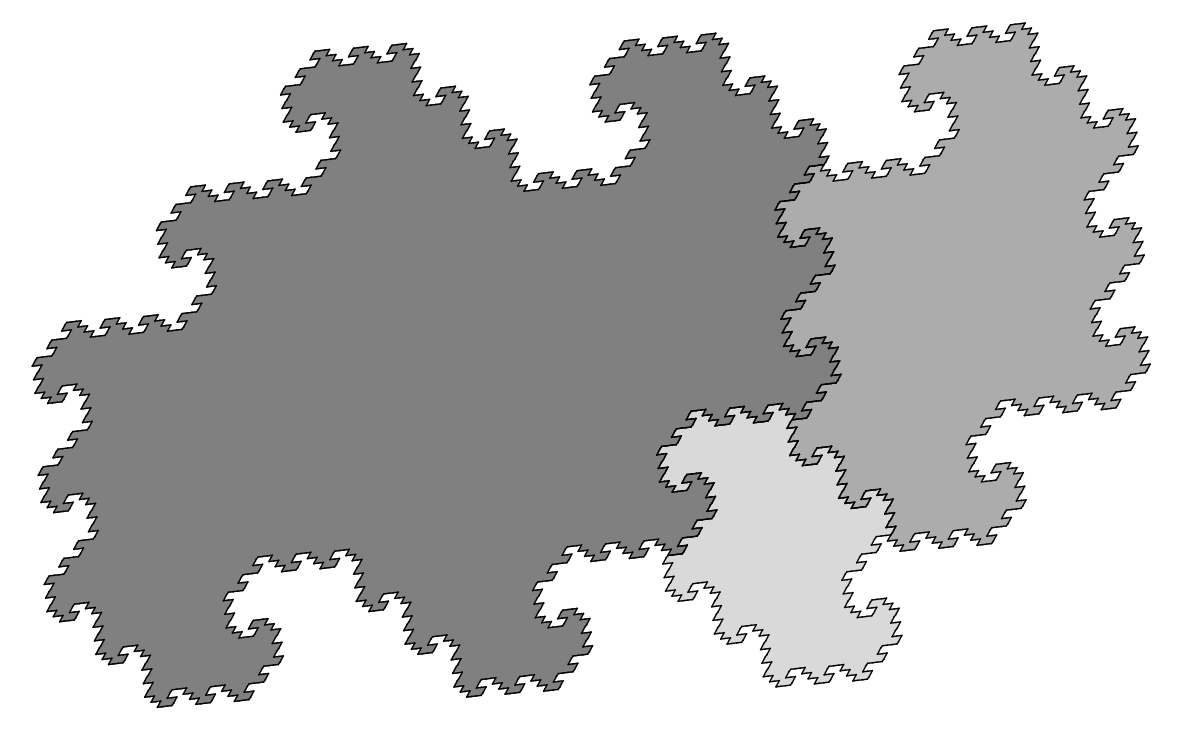}
\hfil
\includegraphics[height=35mm]{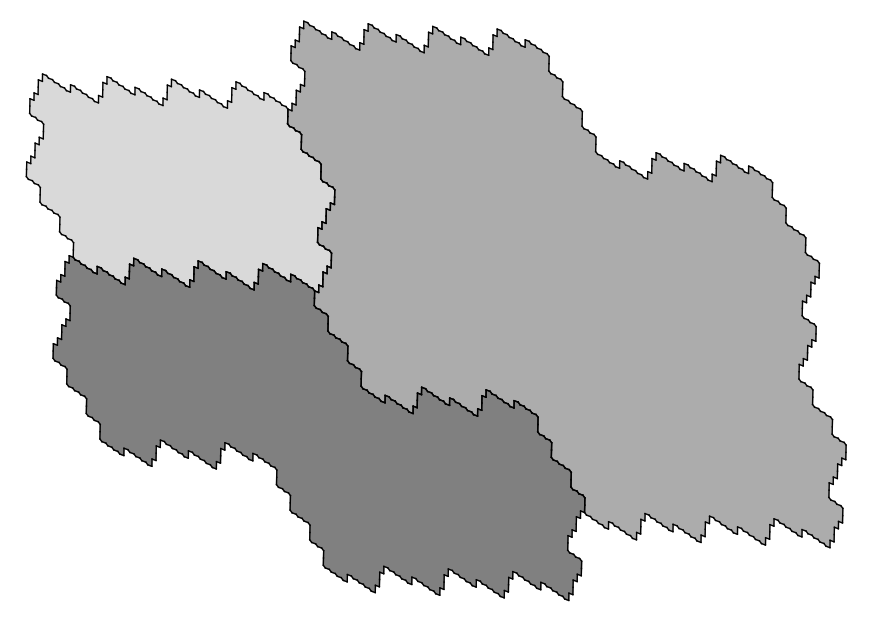}
\caption[]{
Rauzy fractals associated with products
$\sigma_1 \sigma_1 \sigma_2 \sigma_2 \sigma_3 \sigma_3$
and $\sigma_2 \sigma_1 \sigma_3 \sigma_2 \sigma_2 \sigma_2 \sigma_1 \sigma_3$.}
\label{fig:AR}%
\end{figure}

As we will see in the following sections,
the Rauzy fractal associated with such a product of Arnoux-Rauzy substitutions
is always connected (Theorem~\ref{thm:ARfractal}),
but not necessarily simply connected (Section~\ref{sec:counterex}).
Examples of such fractals are depicted in
Figure~\ref{fig:triboexchange}, Figure~\ref{fig:topocool-c} and Figure~\ref{fig:AR}.

In Section~\ref{sec:appl} we will manipulate the dual substitutions $\Sigma_1$, $\Sigma_2$, $\Sigma_3$
associated with the Arnoux-Rauzy substitutions,
which are defined by
$\Sigma_i = \EOS(\sigma_i)$ for $i \in \{1,2,3\}$.
We can compute the $\Sigma_i$ explicitly by using the definition of dual substitutions:
\begin{align*}
\Sigma_1([\bfx, 1]^\star) &= \Msinv \bfx + [\mathbf 0, 1]^\star \cup [\mathbf 0,2]^\star \cup [\mathbf 0, 3]^\star\\
\Sigma_1([\bfx, 2]^\star) &= \Msinv\bfx + [(1,0,0), 2]^\star\\
\Sigma_1([\bfx, 3]^\star) &= \Msinv\bfx + [(1,0,0), 3]^\star
\end{align*}
\begin{align*}
\Sigma_2([\bfx, 1]^\star) &= \Msinv \bfx + [(0,1,0), 1]^\star\\
\Sigma_2([\bfx, 2]^\star) &= \Msinv\bfx +  [\mathbf 0, 1]^\star \cup [\mathbf 0,2]^\star \cup [\mathbf 0, 3]^\star\\
\Sigma_2([\bfx, 3]^\star) &= \Msinv\bfx + [(0,1,0), 3]^\star
\end{align*}
\begin{align*}
\Sigma_3([\bfx, 1]^\star) &= \Msinv \bfx +  [(0,0,1), 1]^\star\\
\Sigma_3([\bfx, 2]^\star) &= \Msinv\bfx + [(0,0,1), 2]^\star\\
\Sigma_3([\bfx, 3]^\star) &= \Msinv\bfx + [\mathbf 0, 1]^\star \cup [\mathbf 0,2]^\star \cup [\mathbf 0, 3]^\star,
\end{align*}
which can be represented graphically as follows,
where the black dot stands for $\bfx$ on the left hand side
and for $\Msinv\bfx$ on the right-hand side.
\[
\definecolor{facecolor}{rgb}{0.8,0.8,0.8}
\renewcommand{\tabcolsep}{1.5mm}
\renewcommand{\arraystretch}{1.2}
\centering
\Sigma_1 : \left\{
    \begin{tabular}{rcl}%
    \myvcenter{%
    \begin{tikzpicture}
    [x={(-0.216506cm,-0.125000cm)}, y={(0.216506cm,-0.125000cm)}, z={(0.000000cm,0.250000cm)}]
    \fill[fill=facecolor, draw=black, shift={(0,0,0)}]
    (0, 0, 0) -- (0, 1, 0) -- (0, 1, 1) -- (0, 0, 1) -- cycle;
    \node[circle,fill=black,draw=black,minimum size=1.2mm,inner sep=0pt] at (0,0,0) {};
    \end{tikzpicture}}%
     & \myvcenter{$\mapsto$} & 
    \myvcenter{%
    \begin{tikzpicture}
    [x={(-0.216506cm,-0.125000cm)}, y={(0.216506cm,-0.125000cm)}, z={(0.000000cm,0.250000cm)}]
    \fill[fill=facecolor, draw=black, shift={(0,0,0)}]
    (0, 0, 0) -- (0, 1, 0) -- (0, 1, 1) -- (0, 0, 1) -- cycle;
    \fill[fill=facecolor, draw=black, shift={(0,0,0)}]
    (0, 0, 0) -- (0, 0, 1) -- (1, 0, 1) -- (1, 0, 0) -- cycle;
    \fill[fill=facecolor, draw=black, shift={(0,0,0)}]
    (0, 0, 0) -- (1, 0, 0) -- (1, 1, 0) -- (0, 1, 0) -- cycle;
    \node[circle,fill=black,draw=black,minimum size=1.2mm,inner sep=0pt] at (0,0,0) {};
    \end{tikzpicture}} \\
    \myvcenter{%
    \begin{tikzpicture}
    [x={(-0.216506cm,-0.125000cm)}, y={(0.216506cm,-0.125000cm)}, z={(0.000000cm,0.250000cm)}]
    \fill[fill=facecolor, draw=black, shift={(0,0,0)}]
    (0, 0, 0) -- (0, 0, 1) -- (1, 0, 1) -- (1, 0, 0) -- cycle;
    \node[circle,fill=black,draw=black,minimum size=1.2mm,inner sep=0pt] at (0,0,0) {};
    \end{tikzpicture}}%
     & \myvcenter{$\mapsto$} & 
    \myvcenter{%
    \begin{tikzpicture}
    [x={(-0.216506cm,-0.125000cm)}, y={(0.216506cm,-0.125000cm)}, z={(0.000000cm,0.250000cm)}]
    \draw[thick, densely dotted] (0,0,0) -- (1,0,0);
    \fill[fill=facecolor, draw=black, shift={(1,0,0)}]
    (0, 0, 0) -- (0, 0, 1) -- (1, 0, 1) -- (1, 0, 0) -- cycle;
    \node[circle,fill=black,draw=black,minimum size=1.2mm,inner sep=0pt] at (0,0,0) {};
    \end{tikzpicture}} \\
    \myvcenter{%
    \begin{tikzpicture}
    [x={(-0.216506cm,-0.125000cm)}, y={(0.216506cm,-0.125000cm)}, z={(0.000000cm,0.250000cm)}]
    \fill[fill=facecolor, draw=black, shift={(0,0,0)}]
    (0, 0, 0) -- (1, 0, 0) -- (1, 1, 0) -- (0, 1, 0) -- cycle;
    \node[circle,fill=black,draw=black,minimum size=1.2mm,inner sep=0pt] at (0,0,0) {};
    \end{tikzpicture}}%
     & \myvcenter{$\mapsto$} & 
    \myvcenter{%
    \begin{tikzpicture}
    [x={(-0.216506cm,-0.125000cm)}, y={(0.216506cm,-0.125000cm)}, z={(0.000000cm,0.250000cm)}]
    \draw[thick, densely dotted] (0,0,0) -- (1,0,0);
    \fill[fill=facecolor, draw=black, shift={(1,0,0)}]
    (0, 0, 0) -- (1, 0, 0) -- (1, 1, 0) -- (0, 1, 0) -- cycle;
    \node[circle,fill=black,draw=black,minimum size=1.2mm,inner sep=0pt] at (0,0,0) {};
    \end{tikzpicture}}
    \end{tabular}
\right.
\quad
\Sigma_2 : \left\{
    \begin{tabular}{rcl}%
    \myvcenter{%
    \begin{tikzpicture}
    [x={(-0.216506cm,-0.125000cm)}, y={(0.216506cm,-0.125000cm)}, z={(0.000000cm,0.250000cm)}]
    \fill[fill=facecolor, draw=black, shift={(0,0,0)}]
    (0, 0, 0) -- (0, 1, 0) -- (0, 1, 1) -- (0, 0, 1) -- cycle;
    \node[circle,fill=black,draw=black,minimum size=1.2mm,inner sep=0pt] at (0,0,0) {};
    \end{tikzpicture}}%
     & \myvcenter{$\mapsto$} & 
    \myvcenter{%
    \begin{tikzpicture}
    [x={(-0.216506cm,-0.125000cm)}, y={(0.216506cm,-0.125000cm)}, z={(0.000000cm,0.250000cm)}]
    \draw[thick, densely dotted] (0,0,0) -- (0,1,0);
    \fill[fill=facecolor, draw=black, shift={(0,1,0)}]
    (0, 0, 0) -- (0, 1, 0) -- (0, 1, 1) -- (0, 0, 1) -- cycle;
    \node[circle,fill=black,draw=black,minimum size=1.2mm,inner sep=0pt] at (0,0,0) {};
    \end{tikzpicture}} \\
    \myvcenter{%
    \begin{tikzpicture}
    [x={(-0.216506cm,-0.125000cm)}, y={(0.216506cm,-0.125000cm)}, z={(0.000000cm,0.250000cm)}]
    \fill[fill=facecolor, draw=black, shift={(0,0,0)}]
    (0, 0, 0) -- (0, 0, 1) -- (1, 0, 1) -- (1, 0, 0) -- cycle;
    \node[circle,fill=black,draw=black,minimum size=1.2mm,inner sep=0pt] at (0,0,0) {};
    \end{tikzpicture}}%
     & \myvcenter{$\mapsto$} & 
    \myvcenter{%
    \begin{tikzpicture}
    [x={(-0.216506cm,-0.125000cm)}, y={(0.216506cm,-0.125000cm)}, z={(0.000000cm,0.250000cm)}]
    \fill[fill=facecolor, draw=black, shift={(0,0,0)}]
    (0, 0, 0) -- (0, 1, 0) -- (0, 1, 1) -- (0, 0, 1) -- cycle;
    \fill[fill=facecolor, draw=black, shift={(0,0,0)}]
    (0, 0, 0) -- (0, 0, 1) -- (1, 0, 1) -- (1, 0, 0) -- cycle;
    \fill[fill=facecolor, draw=black, shift={(0,0,0)}]
    (0, 0, 0) -- (1, 0, 0) -- (1, 1, 0) -- (0, 1, 0) -- cycle;
    \node[circle,fill=black,draw=black,minimum size=1.2mm,inner sep=0pt] at (0,0,0) {};
    \end{tikzpicture}} \\
    \myvcenter{%
    \begin{tikzpicture}
    [x={(-0.216506cm,-0.125000cm)}, y={(0.216506cm,-0.125000cm)}, z={(0.000000cm,0.250000cm)}]
    \fill[fill=facecolor, draw=black, shift={(0,0,0)}]
    (0, 0, 0) -- (1, 0, 0) -- (1, 1, 0) -- (0, 1, 0) -- cycle;
    \node[circle,fill=black,draw=black,minimum size=1.2mm,inner sep=0pt] at (0,0,0) {};
    \end{tikzpicture}}%
     & \myvcenter{$\mapsto$} & 
    \myvcenter{%
    \begin{tikzpicture}
    [x={(-0.216506cm,-0.125000cm)}, y={(0.216506cm,-0.125000cm)}, z={(0.000000cm,0.250000cm)}]
    \draw[thick, densely dotted] (0,0,0) -- (0,1,0);
    \fill[fill=facecolor, draw=black, shift={(0,1,0)}]
    (0, 0, 0) -- (1, 0, 0) -- (1, 1, 0) -- (0, 1, 0) -- cycle;
    \node[circle,fill=black,draw=black,minimum size=1.2mm,inner sep=0pt] at (0,0,0) {};
    \end{tikzpicture}}
    \end{tabular}
\right.
\quad
\Sigma_3 \ : \ \left\{
    \begin{tabular}{rcl}%
    \myvcenter{%
    \begin{tikzpicture}
    [x={(-0.216506cm,-0.125000cm)}, y={(0.216506cm,-0.125000cm)}, z={(0.000000cm,0.250000cm)}]
    \fill[fill=facecolor, draw=black, shift={(0,0,0)}]
    (0, 0, 0) -- (0, 1, 0) -- (0, 1, 1) -- (0, 0, 1) -- cycle;
    \node[circle,fill=black,draw=black,minimum size=1.2mm,inner sep=0pt] at (0,0,0) {};
    \end{tikzpicture}}%
     & \myvcenter{$\mapsto$} & 
    \myvcenter{%
    \begin{tikzpicture}
    [x={(-0.216506cm,-0.125000cm)}, y={(0.216506cm,-0.125000cm)}, z={(0.000000cm,0.250000cm)}]
    \draw[thick, densely dotted] (0,0,0) -- (0,0,1);
    \fill[fill=facecolor, draw=black, shift={(0,0,1)}]
    (0, 0, 0) -- (0, 1, 0) -- (0, 1, 1) -- (0, 0, 1) -- cycle;
    \node[circle,fill=black,draw=black,minimum size=1.2mm,inner sep=0pt] at (0,0,0) {};
    \end{tikzpicture}} \\
    \myvcenter{%
    \begin{tikzpicture}
    [x={(-0.216506cm,-0.125000cm)}, y={(0.216506cm,-0.125000cm)}, z={(0.000000cm,0.250000cm)}]
    \fill[fill=facecolor, draw=black, shift={(0,0,0)}]
    (0, 0, 0) -- (0, 0, 1) -- (1, 0, 1) -- (1, 0, 0) -- cycle;
    \node[circle,fill=black,draw=black,minimum size=1.2mm,inner sep=0pt] at (0,0,0) {};
    \end{tikzpicture}}%
     & \myvcenter{$\mapsto$} & 
    \myvcenter{%
    \begin{tikzpicture}
    [x={(-0.216506cm,-0.125000cm)}, y={(0.216506cm,-0.125000cm)}, z={(0.000000cm,0.250000cm)}]
    \draw[thick, densely dotted] (0,0,0) -- (0,0,1);
    \fill[fill=facecolor, draw=black, shift={(0,0,1)}]
    (0, 0, 0) -- (0, 0, 1) -- (1, 0, 1) -- (1, 0, 0) -- cycle;
    \node[circle,fill=black,draw=black,minimum size=1.2mm,inner sep=0pt] at (0,0,0) {};
    \end{tikzpicture}} \\
    \myvcenter{%
    \begin{tikzpicture}
    [x={(-0.216506cm,-0.125000cm)}, y={(0.216506cm,-0.125000cm)}, z={(0.000000cm,0.250000cm)}]
    \fill[fill=facecolor, draw=black, shift={(0,0,0)}]
    (0, 0, 0) -- (1, 0, 0) -- (1, 1, 0) -- (0, 1, 0) -- cycle;
    \node[circle,fill=black,draw=black,minimum size=1.2mm,inner sep=0pt] at (0,0,0) {};
    \end{tikzpicture}}%
     & \myvcenter{$\mapsto$} & 
    \myvcenter{%
    \begin{tikzpicture}
    [x={(-0.216506cm,-0.125000cm)}, y={(0.216506cm,-0.125000cm)}, z={(0.000000cm,0.250000cm)}]
    \fill[fill=facecolor, draw=black, shift={(0,0,0)}]
    (0, 0, 0) -- (0, 1, 0) -- (0, 1, 1) -- (0, 0, 1) -- cycle;
    \fill[fill=facecolor, draw=black, shift={(0,0,0)}]
    (0, 0, 0) -- (0, 0, 1) -- (1, 0, 1) -- (1, 0, 0) -- cycle;
    \fill[fill=facecolor, draw=black, shift={(0,0,0)}]
    (0, 0, 0) -- (1, 0, 0) -- (1, 1, 0) -- (0, 1, 0) -- cycle;
    \node[circle,fill=black,draw=black,minimum size=1.2mm,inner sep=0pt] at (0,0,0) {};
    \end{tikzpicture}}
    \end{tabular}
\right..
\]

\section{Covering by a set of patterns}
\label{sec:cover}
A \emph{pattern} is a finite set of unit faces.
In the following, we will not make the distinction between a pattern
and the union of its elements.

The aim of this section is to introduce the notion of $\mcL$-\emph{covering}
of a pattern by a set of patterns $\mcL$ (Definition~\ref{def:cover}),
and the notion of \emph{stability} of a set of patterns
with respect to a dual substitution (Definition~\ref{def:stable}).
We then use these concepts to give a simple
sufficient condition for the connectedness of a pattern (Proposition~\ref{prop:connectedim}),
which will be used in Section~\ref{sec:appl}.
The notion of $\mcL$-covering already appeared in~\cite{IO93, IO94, ABI02, ABS04},
but for patterns consisting of two faces only.
We will need to work with larger patterns to be able to deal with Arnoux-Rauzy substitutions;
see also the discussion in Section~\ref{sect:discuss_LAR}.

\begin{defi}[$\mcL$-covering]
\label{def:cover}
Let $D$ be a union of unit faces and $\mcL$ be a set of patterns.
An \emph{$\mcL$-chain} from a face $e \in D$ to a face $f \in D$
is a finite sequence of patterns $(p_1, \ldots, p_n) \in \mcL^n$ such that:
\begin{enumerate}
  \item $e \in p_1$ \quad and \quad $f \in p_n$;
  \item $p_k$ and $p_{k+1}$ have at least one face in common,
    for all $k \in \{1, \ldots, n-1\}$;
  \item $p_k \subseteq D$ \quad for all $k \in \{1, \ldots, n\}$.
\end{enumerate}
We say that $D$ is \emph{$\mcL$-covered} if for all faces $e, f \in D$,
there exists an $\mcL$-chain from $e$ to $f$.
\end{defi}

Roughly speaking, $D$ being $\mcL$-covered means that
we can connect any two faces of $D$ by a ``path'' made of patterns of $\mcL$
in which two consecutive patterns share at least one face.

\begin{figure}[ht]
\centering
\includegraphics{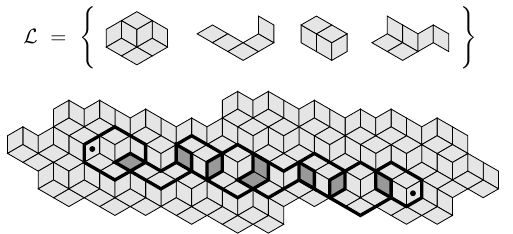}
\caption{Example of an $\mcL$-chain. Dark grey indicates the intersection of two patterns.}
\label{fig:ex_chain}
\end{figure}

The following lemma states that concatenations of $\mcL$-chains remain $\mcL$-chains.

\begin{lem}[\cite{IO94}]
\label{lem:chainconcat}
Let $\mathcal L$ be a set of patterns,
$D$ be a union of unit faces,
and $e, f, g$ three faces of $D$.
If there exists an $\mathcal L$-chain from $e$ to $f$
and an $\mathcal L$-chain from $f$ to $g$,
then there exists an $\mathcal L$-chain from $e$ to $g$.
\end{lem}

\begin{proof}
Let $(p_1, \ldots, p_n)$ be an $\mathcal L$-chain from $e$ to $f$,
and $(q_1, \ldots, q_m)$ an $\mathcal L$-chain from $f$ to $g$.
An $\mathcal L$-chain from $e$ to $g$ is given by
$(p_1, \ldots, p_n, q_1, \ldots, q_m)$, because $f$ is in $p_n \cap q_1$.
\end{proof}

Let $\Sigma$ be a dual substitution and $D$ be an $\mcL$-covered pattern.
We want to know when $\Sigma(D)$ is $\mcL$-covered.
It turns out that there is a simple sufficient condition to check this,
namely the \emph{stability} of $\mcL$ under $\Sigma$ (Definition~\ref{def:stable}),
as stated in Proposition~\ref{prop:coverprop} below.

\begin{defi}[Stability]
\label{def:stable}
Let $\Sigma$ be a dual substitution.
A set of patterns $\mathcal L$ is \emph{stable} under $\Sigma$
if $\Sigma(p)$ is $\mathcal L$-covered for all $p \in \mathcal L$.
\end{defi}

\begin{prop}
\label{prop:coverprop}
Let $\mathcal L$ be a set of patterns
that is stable under a dual substitution $\Sigma$,
and let $D$ be an $\mathcal L$-covered union of unit faces.
Then $\Sigma(D)$ is $\mathcal L$-covered.
\end{prop}

\begin{figure}[ht]
\centering
\includegraphics{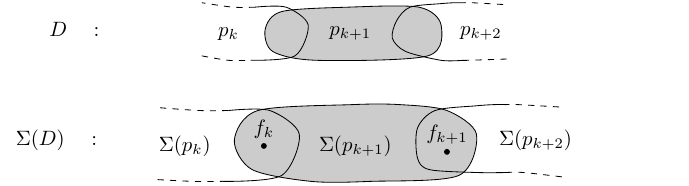}
\caption{Illustration of the proof of Proposition~\ref{prop:coverprop}.}
\label{fig:coverprop}
\end{figure}

\begin{proof}
Let $f$ and $f'$ be two faces of $\Sigma(D)$.
To prove that $\Sigma(D)$ is $\mathcal L$-covered,
we need to construct an $\mathcal L$-chain from $f$ to $f'$.
Let $e$ and $e'$ be two faces of $D$ such that $f \in \Sigma(e)$ and $f' \in \Sigma(e')$.
Since $D$ is $\mathcal L$-covered,
there exists an $\mathcal L$-chain $(p_1, \ldots, p_n)$ from $e$ to $e'$.
For all $k \in \{2, \ldots, n-1\}$,
let $f_k$ be a face of $\Sigma(p_k \cap p_{k+1})$,
and let $f_1 = f$, $f_n = f'$.
Such a face $f_k$ exists because $(p_1, \ldots, p_n)$ is an $\mathcal L$-chain.

For all $k \in \{1, \ldots, n-1\}$,
there exists an $\mathcal L$-chain from $f_k$ to $f_{k+1}$,
because $f_k$ and $f_{k+1}$ are in $\Sigma(p_{k+1})$
and $\mathcal L$ is stable under $\Sigma$ (see Figure~\ref{fig:coverprop}).
Lemma~\ref{lem:chainconcat} implies that
the concatenation of the $\mathcal L$-chains from $f_k$ to $f_{k+1}$
is an $\mathcal L$-chain from $f$ to $f'$.
\end{proof}

The stability of a given set $\mcL$ under $\Sigma$ is easy to verify:
it can be done in finite time, since there are only a finite number of patterns to check,
and a pattern can be covered only in finitely many ways by patterns of $\mcL$.

The following proposition relates $\mcL$-coverings and connectedness:
if a pattern $D$ is $\mcL$-covered by a set of patterns $\mcL$
whose elements have a connected projection by $\pic$,
then $\pic(D)$ is also connected.

\begin{prop}
\label{prop:connectedim}
Let $\mathcal L$ be a set of patterns
and $D$ be an $\mathcal L$-covered union of unit faces.
If $\pic(p)$ is connected for all $p \in \mathcal L$,
then $\pic(D)$ is connected.
\end{prop}

\begin{proof}
Let $x, y$ be two points of $\pic(D)$ and $e, f$ two faces of $D$
such that $x \in \pic(e)$ and $y \in \pic(f)$.
Let $(p_1, \ldots, p_n)$ be an $\mathcal L$-chain from $e$ to $f$.
The sets $\pic(p_i)$ are connected and $\pic(p_i) \cup \pic(p_{i+1}) \neq \varnothing$
for all $i$, so there exists a path from $x \in \pic(p_1)$ to $y \in \pic(p_n)$.
\end{proof}

The following basic lemma will be useful in the next section.

\begin{lem}
\label{lem:connectedlimit}
Let $K_1, K_2, \ldots$ be a sequence of compact subsets of $\bbR^2$
that converges to a Hausdorff limit $K$.
If the $K_n$ are connected, then $K$ is connected.
\end{lem}

\begin{proof}
Suppose that $K$ is not connected:
let $A$ and $B$ be two disjoint non-empty closed sets such that $K = A \cup B$.
Let $\varepsilon = \min\{\|x-y\| : (x,y) \in A \times B\}$,
and let $n \geq 0$ such that $\dH(K_n, K) \leq \varepsilon / 3$.
Then we have $K_n = (\mathcal V_{\varepsilon/3}(A) \cap K_n) \cup (\mathcal V_{\varepsilon/3}(B) \cap K_n)$,
where $\mathcal V_\varepsilon(X)$ stands for the $\varepsilon$-neighborhood of a set $X$.
Hence $K_n$ is not connected because it is the union of two disjoint non-empty closed sets:
a contradiction.
\end{proof}

\section{Applications to Arnoux-Rauzy substitutions}
\label{sec:appl}
\subsection{Main results}
In this section, we use the tools developed in Section~\ref{sec:cover}
to prove the connectedness of the images of
$\mcU = [{\bf 0},1]^\star \cup [{\bf 0},2]^\star \cup [{\bf 0},3]^\star$
under any finite product of the dual substitutions $\Sigma_1$, $\Sigma_2$, $\Sigma_3$ defined by
$\Sigma_i = \EOS(\sigma_i)$ for $i \in \{1,2,3\}$.

To this end, we introduce a finite set of connected patterns
$\LAR$ (Equation~\ref{eq:LAR}) which covers $\mcU$,
and which is stable under $\Sigma_1$, $\Sigma_2$ and $\Sigma_3$ (Proposition~\ref{prop:stableAR}).
The covering property will then be transferred to all the forward images of $\mcU$,
which yields their connectedness (Theorem~\ref{thm:ARapprox}).
Let
\begin{eqnarray}
\label{eq:LAR}
\LAR \quad = \quad
\left\{\raisebox{7mm}{}\right.
\myvcenter{\raisebox{2mm}{\includegraphics[scale=0.9]{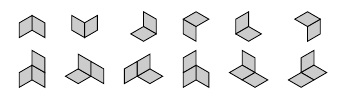}}\hspace{-5pt}}
\left.\raisebox{7mm}{}\right\}.
\end{eqnarray}
We do not explicit define each pattern of $\LAR$,
since it can clearly be deduced from the above graphical representation
(we require that each pattern is a connected subset of $\bbR^3$).
We also do not specify the position in $\bbZ^3$ of each pattern,
because it does not matter:
any choice is compatible with the proofs below.

\begin{prop}
\label{prop:stableAR}
The set of patterns $\LAR$ is stable under
$\Sigma_1$, $\Sigma_2$ and $\Sigma_3$.
\end{prop}

\begin{proof}
We must prove that every pattern of
$\Sigma_i(\LAR)$ is covered by $\LAR$, for $i = 1,2,3$,
which makes a total of $36$ patterns to check:
\begin{align*}
\Sigma_1(\LAR) &=
  \left\{\raisebox{7mm}{}\right.
  \myvcenter{\raisebox{1.5mm}{\includegraphics[scale=0.8]{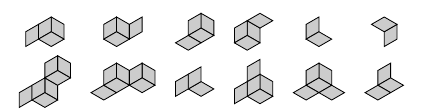}}\hspace{-5pt}}
  \left.\raisebox{7mm}{}\right\}; \displaybreak[0] \\
\Sigma_2(\LAR) &=
  \left\{\raisebox{7mm}{}\right.
  \myvcenter{\raisebox{1.5mm}{\includegraphics[scale=0.8]{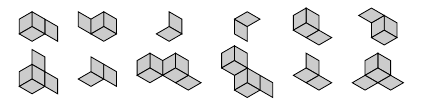}}\hspace{-5pt}}
  \left.\raisebox{7mm}{}\right\}; \displaybreak[0] \\
\Sigma_3(\LAR) &=
  \left\{\raisebox{7mm}{}\right.
  \myvcenter{\raisebox{1.5mm}{\includegraphics[scale=0.8]{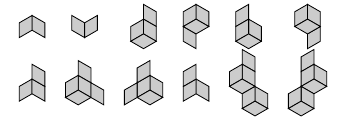}}\hspace{-5pt}}
  \left.\raisebox{7mm}{}\right\}.
\end{align*}
This can easily be checked for each of these patterns,
as for the following pattern of $\Sigma_1(\LAR)$, for example:
\[
\renewcommand{\tabcolsep}{0.3mm}
\begin{tabular}{ccccccc}
\myvcenter{\includegraphics{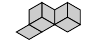}} & \enspace : \quad &
\myvcenter{\includegraphics{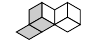}} &
\myvcenter{\includegraphics{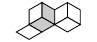}} &
\myvcenter{\includegraphics{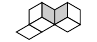}} &
\myvcenter{\includegraphics{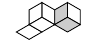}} &
\myvcenter{\includegraphics{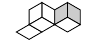}}
\end{tabular}.
\]
From this graphical representation,
we can deduce that there is an $\LAR$-chain between each two faces
of the pattern.
All the other patterns of $\Sigma_1(\LAR)$ are also $\LAR$-covered,
and the patterns of the sets $\Sigma_2(\LAR)$ and $\Sigma_3(\LAR)$
(which are symmetrical copies of the patterns of $\Sigma_1(\LAR)$)
admit similar $\LAR$-coverings.
\end{proof}

We now use the stability of $\LAR$ and the results of Section~\ref{sec:cover}
to obtain our connectedness results.

\begin{thm}
\label{thm:ARapprox}
The set $\pic(\Sigma_{i_1} \cdots \Sigma_{i_n}(\mcU))$ is connected
for every $i_1, \ldots, i_n \in \{1, 2, 3\}$.
\end{thm}

\begin{proof}
Proposition~\ref{prop:coverprop} implies that
$\Sigma_{i_1} \cdots \Sigma_{i_n}(\mcU)$ is $\LAR$-covered,
because $\LAR$ is stable under the $\Sigma_i$ (Proposition~\ref{prop:stableAR})
and $\mcU$ is $\LAR$-covered.
The projection by $\pic$ of every pattern of $\LAR$ is connected,
so Proposition~\ref{prop:connectedim} implies that
the set $\pic(\Sigma_{i_1} \cdots \Sigma_{i_n}(\mcU))$ is connected.
\end{proof}

Since connectedness is preserved in the Hausdorff limit (Lemma~\ref{lem:connectedlimit}),
it follows that the fractals associated with Arnoux-Rauzy substitutions are connected.

\begin{thm}
\label{thm:ARfractal}
Let $\sigma$ be a Pisot irreducible finite product of $3$-letter Arnoux-Rauzy substitutions.
Then, the Rauzy fractal of $\sigma$ is connected.
\end{thm}

As a last application we obtain the following theorem,
stating properties about Markov partitions for the toral automorphism
defined by the action of $\bfM_\sigma$ on $\mathbb T^3$,
where $\sigma$ is a product of Arnoux-Rauzy substitutions.
(See~\cite{LM95} for some background on Markov partitions.)

The existence of such Markov partitions is due to Bowen~\cite{Bow08},
but this construction is not explicit,
and the boundaries of the partitions are necessarily of fractal nature~\cite{Bow78} in three dimensions or more.
In the specific two-dimensional case,
explicit partitions can be constructed~\cite{Adl98}:
two rectangles are always enough.

By combining the results from~\cite{IO93,Pra99,Sie00}
(which establish links between the Rauzy fractal of $\sigma$ and Markov partitions of $\bfM_\sigma$),
and the results from~\cite{BJS12,BSW13}
(which establish the pure discrete spectrum of Arnoux-Rauzy substitutions),
we can apply Theorem~\ref{thm:ARfractal} to obtain the following result.

\begin{thm}
\label{thm:ARmarkov}
Let $\sigma$ be a Pisot irreducible finite product of $3$-letter Arnoux-Rauzy substitutions.
Then, there exists a partition $(\mathcal P_1, \mathcal P_2, \mathcal P_3)$
of the $3$-dimensional torus $\mathbb T^3$
which is a Markov partition for the toral automorphism $(\mathbb T^3, \bfM_\sigma)$
provided by the incidence matrix $\Ms$ of $\sigma$,
and such that the elements of the partitition $\mathcal P_1, \mathcal P_2, \mathcal P_3$ are connected.
\end{thm}

\subsection{About the set $\LAR$}
\label{sect:discuss_LAR}
Our main result (Theorem~\ref{thm:ARapprox}),
which states the stability of $\LAR$ by $\Sigma_1$, $\Sigma_2$ and $\Sigma_3$,
is based on a precise choice of the patterns in $\LAR$.

A possible approach to find such a suitable set of patterns $\LAR$
is to start with all the two-face edge-connected patterns,
compute their images by the $\Sigma_i$,
and detect which images are \emph{not} edge-connected, such as:
\[
\includegraphics{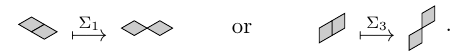}
\]
In the present case, it is possible to ``complete''
each such pattern, in such a way that its image is edge-connected:
\[
\includegraphics{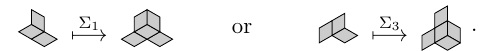}
\]
The resulting patterns not only have edge-connected images by the $\Sigma_i$,
but also have the property that their images are all $\LAR$-covered.

\section{Further questions and counterexamples}
\label{sec:counterex}

In this section, we provide examples to show that
some possible generalizations or extensions of Theorem~\ref{thm:ARfractal} are not true.

\subsection{The converse of Theorem~\ref{thm:ARfractal}}
We give an example of a substitution $\sigma$ whose Rauzy fractal is connected
and such that no power of $\sigma$ can be written as a product of Arnoux-Rauzy substitutions.
This proves that the converse of Theorem~\ref{thm:ARfractal} does not hold, \emph{i.e.},
it is not true that every connected Rauzy fractal is
associated with a product of Arnoux-Rauzy substitutions.

\begin{figure}[ht]%
\centering
\includegraphics[height=30mm]{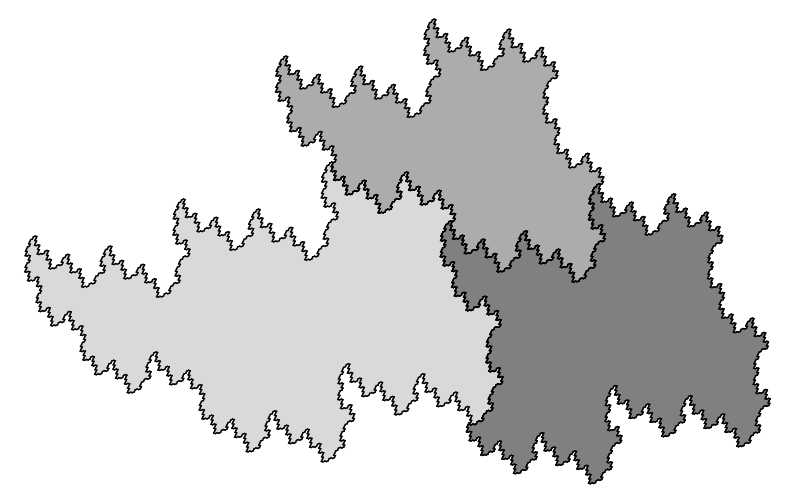}
\hfil
\includegraphics[height=30mm]{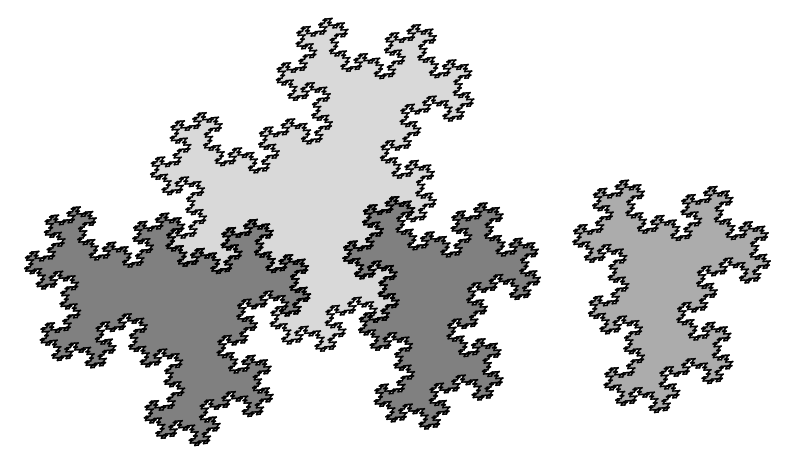}
\caption[]{
Left:
a simply connected Rauzy fractal
where no power of $\sigma$ can be written as a product of Arnoux-Rauzy substitutions
($1 \mapsto 32131, 2 \mapsto 321, 3 \mapsto 3213$).
Right:
a disconnected Rauzy fractal
where $\sigma$ is a product of elementary substitutions
($1 \mapsto 31, 2 \mapsto 12, 3 \mapsto 31123$).
}
\label{fig:counterex}%
\end{figure}

Let $\sigma$ be the primitive substitution defined by $1 \mapsto 32131$, $2 \mapsto 321$, $3 \mapsto 3213$.
It is easy to check that there are at least $14$ different words of length $6$
in an infinite fixed point of $\sigma$.
Hence, no power of $\sigma$ is a product of Arnoux-Rauzy substitutions,
because Arnoux-Rauzy sequences have factor complexity $2n+1$.
It can be checked algorithmically that the Rauzy fractal of $\sigma$ is connected,
using the methods described in~\cite{ST10}.
Let us remark that these methods also enable us to prove that
the Rauzy fractal of $\sigma$ is also \emph{simply} connected,
which makes the counterexample even stronger.
The corresponding Rauzy fractal is shown in Figure~\ref{fig:counterex}.

Many other examples can easily be found,
such as the substitution given in Figure~\ref{fig:topocool-a},
whose Rauzy fractal is connected (but not simply connected).

\subsection{Simple connectedness and Arnoux-Rauzy substitutions}
Simple connectedness of the Rauzy fractal does not hold in general
for products of Arnoux-Rauzy substitutions.
Indeed, the fractal associated with the substitution
$\sigma \ = \ \sigma_1^4 \sigma_2^4 \sigma_3^4$
is not simply connected because $\sigma$ is equal to
the cube of the substitution given in Figure~\ref{fig:topocool-c},
so it has the same Rauzy fractal,
which has uncountable fundamental group~\cite{ST10}.

\subsection{Connectedness and invertible substitutions}
In the case of a $2$-letter substitution $\sigma$,
it has been shown that the Rauzy fractal of $\sigma$ is connected if and only if
$\sigma$ is invertible~\cite{EI98}.
(A substitution is \emph{invertible} if it extends to an automorphism of the free group.)

We can see Theorem~\ref{thm:ARfractal} as a partial analogue of this fact for $3$-letter substitutions,
because Arnoux-Rauzy substitutions are invertible.
The following example shows that Theorem~\ref{thm:ARfractal}
does not hold for a larger class of invertible substitutions,
namely the products of \emph{elementary} substitutions defined by
\[
\varepsilon_{i,j} \ : \ \left \{
  \begin{array}{l}
    j \mapsto ij \\
    k \mapsto k \text{ if } k \neq j
  \end{array}
\right.
\qquad \qquad
\]
for $i, j \in \{1,2,3\}$ and $i \neq j$.
Indeed, every finite product of $\varepsilon_{i,j}$ is invertible,
and we have
$\sigma_1 = \varepsilon_{1,2} \varepsilon_{1,3}$,
$\sigma_2 = \varepsilon_{2,1} \varepsilon_{2,3}$, and
$\sigma_3 = \varepsilon_{3,1} \varepsilon_{3,2}$.
But the Rauzy fractal associated with the substitution
\[
\varepsilon_{1,2} \varepsilon_{3,1} \varepsilon_{2,3} \varepsilon_{1,3} \ : \ \left \{
  \begin{array}{l}
    1 \mapsto 13 \\
    2 \mapsto 21 \\
    3 \mapsto 32113
  \end{array}
\right.
\]
is not connected, as can be checked using the algorithms given in~\cite{ST10}.
This fractal is depicted in Figure~\ref{fig:counterex}.
Let us mention that invertible substitutions on three letters have been characterized~\cite{TWZ04},
and that the finite products of $\varepsilon_{i,j}$ constitute only a proper subclass of invertible substitutions.

\section{Conclusion}
We have given a combinatorial proof of the connectedness of the Rauzy fractals
associated with finite products of Arnoux-Rauzy substitutions.
To do so, we have extended combinatorial techniques from \cite{IO94} ($\mcL$-coverings)
by considering coverings by patterns of more than two faces.

There is a lot of room for future work.
We believe that the techniques introduced in this article will
allow us to prove topological properties
for some other families of Rauzy fractals.
It would also be interesting to characterize the products of
Arnoux-Rauzy substitutions which have simply connected Rauzy fractal.

There are also many interesting related decidability questions:
given a unimodular Pisot irreducible substitution,
is its Rauzy fractal connected? Simply connected?
Is the origin an inner point?
Does it verify the tiling property?
Some of these questions have been addressed (see~\cite{ST10}),
but the techniques used rely on incidence graphs of the subtiles of the fractals.
We would like to investigate and revisit these questions
using the techniques developed in this article.

\bibliographystyle{amsalpha}
\bibliography{biblio}
\label{sec:biblio}
\end{document}